\documentclass{amsart}

\usepackage{a4wide}
\usepackage{graphicx}
\usepackage{amstext, amsmath, amssymb}
\usepackage{booktabs}
\usepackage{url}
\usepackage{import}
\usepackage{amsfonts}
\usepackage{amsbsy}
\usepackage{fancyvrb}
\usepackage{stmaryrd}
\usepackage{pdfsync}
\usepackage{graphicx}
\usepackage{algorithm}
\usepackage[shadow]{todonotes}
\usepackage[numbers]{natbib}

\numberwithin{equation}{section}
\numberwithin{figure}{section}
\numberwithin{table}{section}

\numberwithin{theorem}{section}

\newcommand{\R}{\mathbb{R}}

\newcommand{\foralls}{\forall\,}

\DeclareMathOperator{\tr}{tr}
\DeclareMathOperator{\Ext}{Ext}

\newcommand{\jump}[1]{[#1]}

\DefineVerbatimEnvironment{code}{Verbatim}{frame=single,rulecolor=\color{blue}}

\newcommand{\OS}{\Omega^{s}}     
\newcommand{\OF}{\Omega^{f}}     
\newcommand{\IFF}{\Gamma^{{ff}}}  
\newcommand{\IFS}{\Gamma^{{fs}}}  
\newcommand{\mesh}{\mathcal{T}}
\newcommand{\meshO}{\mesh_0}
\newcommand{\reducedmesh}{\mesh_1^{\ast}}
\newcommand{\reduceddomain}{\Omega_1^{\ast}}
\newcommand{\fluidmesh}{\mesh^f_2}
\newcommand{\solidmesh}{\mesh^s}
\newcommand{\fluidsolidmesh}{\mesh^{fs}}
\newcommand{\fluidsoliddomain}{\Omega^{fs}}
\newcommand{\fluidoverlapregion}{\Omega_O}

\newcommand{\bfu}{\boldsymbol{u}}
\newcommand{\bff}{\boldsymbol{f}}

\newcommand{\bfg}{\boldsymbol{g}}
\newcommand{\bft}{\boldsymbol{t}}
\newcommand{\bfn}{\boldsymbol{n}}
\newcommand{\bfv}{\boldsymbol{v}}
\newcommand{\bfx}{\boldsymbol{x}}
\newcommand{\bfF}{\boldsymbol{F}}
\newcommand{\bfE}{\boldsymbol{E}}
\newcommand{\bfzero}{\boldsymbol{0}}
\newcommand{\bphi}{\boldsymbol{\phi}}

\newcommand{\uF}{\bfu^f}      
\newcommand{\pF}{p^{f}}        
\newcommand{\uS}{\bfu^s}      
\newcommand{\bfsigma}{\boldsymbol \sigma}
\newcommand{\sigmaF}{\bfsigma^f} 
\newcommand{\sigmaS}{\bfsigma^s} 

\newcommand{\refe}[1]{\widehat{#1}}

\newcommand{\meshFk}{\mesh^{f,k+1}}

\newcommand{\meshFtwok}{\mesh^{f,k+1}_2}

\newcommand{\tn}{|\mspace{-1mu}|\mspace{-1mu}|} 
\newcommand{\nablan}{\partial_{\bfn}}             
\newcommand{\meanvalue}[1]{\langle #1 \rangle}

\newcommand{\tab}{\hspace*{2em}}

\begin{document}
\title[\bf A cut and composite mesh method for fluid--structure
interaction] {\bf A Nitsche-based cut finite element method for a
  fluid--structure interaction problem}

\author{Andr\'e Massing \and Mats G.\ Larson \and Anders Logg \and
  Marie E.\ Rognes}

\date{Received: \today / Accepted: }

\begin{abstract}
  We present a new composite mesh finite element method for
  fluid--structure interaction problems. The method is based on
  surrounding the structure by a boundary-fitted fluid mesh which is
  embedded into a fixed background fluid mesh. The embedding allows
  for an arbitrary overlap of the fluid meshes. The coupling between
  the embedded and background fluid meshes is enforced using a
  stabilized Nitsche formulation which allows us to establish
  stability and optimal order \emph{a priori} error estimates,
  see~\cite{MassingLarsonLoggEtAl2013}. We consider here a steady state
  fluid--structure interaction problem where a hyperelastic structure
  interacts with a viscous fluid modeled by the Stokes equations. We
  evaluate an iterative solution procedure based on splitting and
  present three-dimensional numerical examples.  
\end{abstract}

\maketitle

\section{Introduction}

In fluid--structure interaction applications, the underlying geometry
of the computational domain may change significantly due to
displacement of the structure. In order to deal with this situation in
a standard setting with conforming elements, a mesh motion algorithm
must be used. If the displacements are significant, the deformation of
the mesh may lead to deteriorating mesh quality which may ultimately
require re-meshing of the computational domain. Alternative, more
flexible, techniques are therefore of significant practical
interest.

In this paper, we consider a combination of standard moving meshes and
so-called CutFEM
technology~\cite{BurmanClausHansboEtAl2014}. Essentially, the
structure or elastic solid is first embedded into a boundary-fitted
fluid mesh which moves along with the deformation of the solid to keep
the fluid--structure interface intact. The motion of the fluid mesh
surrounding the structure is obtained by solving an elasticity problem
with given displacement at the fluid--structure interface.  The
boundary-fitted fluid mesh is then embedded into a fixed background
mesh where we allow for an arbitrary overlap of the fluid meshes in
order to facilitate the repositioning of the moving fluid mesh within
the fixed background mesh. The fluid is then discretized on both the
moving overlapping domain, using an Arbitrary-Lagrange-Eulerian (ALE)
type
approach\cite{DoneaGiulianiHalleux1982,DoneaHuertaPonthotEtAl2004},
and on the fixed background mesh, using a standard discretization
posed in an Eulerian frame.

The coupling at the fluid--fluid interface between the overlapping and
underlying fluid meshes is handled using a stabilized Nitsche method
developed for the Stokes problem in~\cite{MassingLarsonLoggEtAl2013}. The
stabilization is constructed in such a way that the resulting scheme
is inf-sup stable and the resulting stiffness matrix is
well-conditioned independent of the position of the overlapping fluid
mesh relative to the fixed background fluid mesh. As a result, optimal
order error estimates are also established. In order to deal with the
cut elements arising at the interface, we compute the polyhedra
resulting from the intersections between the overlapping and
background meshes. These polyhedra may then be described using a
partition into tetrahedra; this partition may in turn be used to
perform numerical quadrature. We refer to~\cite{Massing2012a} for a
detailed discussion of the implementation aspects of cut element
techniques in three spatial dimensions. We remark that Nitsche-based
formulations for Stokes boundary and interface problems where the
surface in question is described independently of a single, fixed
background mesh were proposed
in~\cite{BurmanHansbo2013,MassingLarsonLoggEtAl2013a,HansboLarsonZahedi2014a,BurmanClausMassing2014a}.
A Nitsche-based composite mesh method
was first introduced for elliptic problems
in~\cite{HansboHansboLarson2003}.

One may also consider formulations, where the structure is described
via its moving boundary which is immersed into a fixed background
fluid mesh. Prominent examples are Cartesian grid methods,
e.g.~\cite{MurmanAftosmisBerger2003}, the classical immersed boundary
method introduced by~\citet{Peskin1977,Peskin2003}, its finite element
pendant proposed in
\citep{BoffiGastaldi2003,ZhangGerstenbergerWangEtAl2004,ZhangGay2007},
formulations based on Lagrange multipliers,
cf.~\citep{ZhangGerstenbergerWangEtAl2004,Yu2005,GerstenbergerWall2008,GerstenbergerWall2008a,PusoKokkoSettgastEtAl2014}
and on Nitsche's method~\cite{HansboHermansson2003}. However, the use
of an additional boundary-fitted fluid mesh as in the current work is
attractive since it allows for the resolution of boundary layers and
computation of accurate boundary stresses. Often, the construction of
the surrounding fluid mesh can easily be generated by extending the
boundary mesh in the normal direction. We plan to further investigate
the properties of the fluid--structure coupling in future work.

As our proposed scheme combines an ALE-based discretization on the
fluid mesh surrounding the structure with an Eulerian-based
discretization on the fixed background fluid mesh, it can be
classified as hybrid Eulerian-ALE or Chimera approach. Such hybrid
schemes are built upon the concept of overlapping meshes introduced
for finite differences and finite volume schemes in the early works
of~\citet{Volkov1968}, \citet{Starius1977,Starius1980},
\citet{Steger1983} and later by~\citet{ChesshireHenshaw1990}
and~\citet{AftosmisBergerMelton1998}, where the primary concern was to
ease the burden of mesh generation by composing individually meshed,
static geometries. The idea of gluing meshes together was then
explored for finite element methods
by~\citet{CebralLoehner2005,LoehnerCebralCamelliEtAl2007,LoehnerCebralCamelliEtAl2008}
to study the flow around independently meshed complex objects such as
cars, collection of buildings or stents in aortic vessels.  In these
works, relatively simple interpolation schemes were used to
communicate the solution between overlapping meshes. To achieve a
physically more consistent coupling between the solution parts
presented on different domains, Schwarz-type domain iteration schemes
using Dirichlet/Neumann and Robin coupling on overlapping domains have
been proposed for the Navier-Stokes equations
in~\cite{HouzeauxCodina2003}. A completely different route was taken
by~\citet{DayBochev2008} who reformulated elliptic interface problems
as suitable first-order systems augmented with least-square
stabilizations to enforce the interface conditions between the mesh
domains to be tied together.

Introducing special interpolation stencils close to the fluid--fluid
interface, a finite volume based Chimera method for flow problems
involving multiple moving rigid bodies was formulated
in~\cite{Wang2000,EnglishQiuYuEtAl2013}
and~\cite{HenshawSchwendeman2006}, where higher-order Godunov fluxes
where used. This method was then extended
by~\citet{BanksHenshawSchwendeman2012} to deal with (linearly) elastic
solids in two space dimensions, and thus represents an instance of a
hybrid ALE-fixed grid method. This approach has barely been explored
in the context of finite-element methods for fluid--structure
interaction problems: \citet{WallGerstenbergerGamnitzerEtAl2006} and
later \citet{Shahmiri2011} used interpolation between fluid meshes and
extended finite element techniques to couple fluid--fluid meshes,
\citet{BaigesCodina2009} introduced an auxiliary ALE step to convect
information on the fixed background mesh between two consecutive
time-steps.

In contrast to these contributions, our method is based on a
variational finite element approach that leads to a monolithic and
physically consistent coupling between the overlapping and underlying
fluid meshes, which eliminates the need of introducing inconsistent
interpolation operators. In addition, opposed to similar finite
element based approaches presented
e.g. in~\cite{WallGerstenbergerGamnitzerEtAl2006,Shahmiri2011}, our
scheme used for the fluid problem is proven stable and optimally
convergent, even for higher-order elements, independent of the
location of the interface as shown in~\cite{MassingLarsonLoggEtAl2013}.
Thus, the new scheme for the fluid--structure interaction problem
proposed in this work exhibits the necessary robustness that is
essential for developing reliable hybrid ALE-fixed mesh methods.

In the current work, we consider the steady state deformation of a
hyperelastic solid immersed into a viscous fluid governed by the
Stokes equations. We solve for the steady state solution using a fixed
point iteration where in each iteration the fluid, solid, and mesh
motion problems are solved sequentially. We present two numerical
examples in three dimensions, including one example with a
manufactured reference solution.

The outline of the remainder of this paper is as follows: in Section
2, we summarize the governing equations of the fluid--structure
interaction (FSI) problem; in Section 3, we describe the overlapping
mesh method; in Section 4, we present an algorithm for the solution of
the stationary fluid--structure interaction model problem; in Section
5, we present three-dimensional numerical examples; before drawing
some conclusions in Section 6.

\section{A stationary fluid--structure interaction problem}

\begin{figure}
  \begin{center}
    \includegraphics[width=0.45\textwidth]{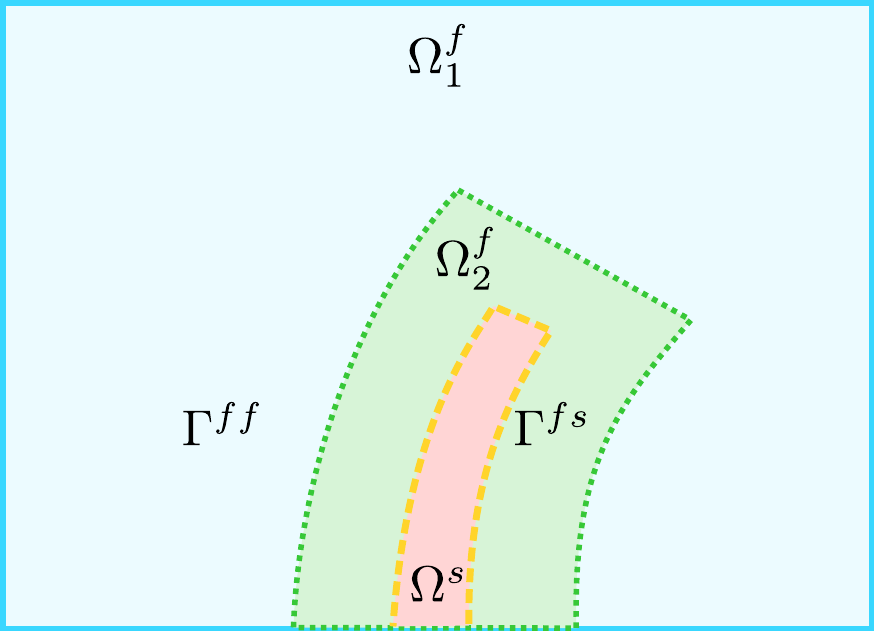}
    \caption{Fluid and structure domains for the stationary
      fluid--structure interaction problem.}
    \label{fig:ref-phy-domains}
  \end{center}
\end{figure}

We consider a fluid--structure interaction problem posed on a domain
$\Omega = \OF \cup \OS$ where $\OF$ is the domain occupied by the
fluid and $\OS$ is the domain occupied by the solid. We assume that
both $\OF$ and $\OS$ are open and bounded and that they are such that
$\OF \cap \OS = \emptyset$. Furthermore, we decompose the fluid domain
into two disjoint subdomains $\OF_1$, $\OF_2$ such that $\OF = \OF_1
\cup \OF_2$. Here, $\OF_2$ represents a part of the fluid domain
surrounding the solid domain $\OS$; more precisely, we assume that
$\partial \OF_1 \cap \partial \OS = \emptyset$. The fluid--structure
interface is denoted by $\IFS = \partial \OF_2 \cap \partial \OS$ and
the interface between the two fluid domains is denoted by $\IFF =
\partial \OF_1 \cap \partial \OF_2$.
Here, the topological boundary $\partial X$ for any given set $X$ is defined by $\partial X = \overline{X} \setminus
\overset{\circ}{X}$
 where $\overline{X}$ and
$\overset{\circ}{X}$ denotes the closure and interior of $X$, respectively.
 For simplicity, we assume that
the fluid domain boundary consists of two disjoint parts: $\partial
\OF = \IFS \cup \partial \OF_D$, and that the solid domain boundary
decomposes in a similar manner: $\partial \OS = \IFS \cup \partial
\OS_D$. This notation is summarized in
Figure~\ref{fig:ref-phy-domains}.

We assume that the fluid dynamics are governed by the Stokes equations
of the following form: find the fluid velocity $\uF : \OF \to \R^3$
and the fluid pressure $\pF: \OF \to \R$ such that
\begin{alignat}{3}
  - \nabla \cdot (\nu^f \nabla \uF - \pF \mathbf{I}) &= \bff^f & &
  \quad \text{in } \OF,
  \label{eq:fluid-momentum-pde}
  \\
  \nabla \cdot \uF &= 0 & &\quad \text{in } \OF,
  \label{eq:fluid-compressibility-pde}
\end{alignat}
where $\bff^f$ is a given body force and $\nu^f$ is the fluid
viscosity.

Next, we assume that the velocity is prescribed on both the
fluid--structure interface and on the remainder of the fluid boundary:
\begin{alignat}{3}
  \uF &= 0 && \quad \text{on } \IFS ,
  \label{eq:fluid-bc1}
  \\
  \uF &= \bfg^f && \quad \text{on } \partial \OF_D .
  \label{eq:fluid-bc2}
\end{alignat}
Moreover, we enforce the continuity of the fluid velocity and of the
fluid ``stress'' on the fluid--fluid interface by the following
conditions:
\begin{alignat}{3}
  \jump{\uF} &= 0 & & \quad \text{on } \IFF,
  \label{eq:fluid-interface-jump-1}
  \\
  \jump{ \left (\nu^f \nabla \uF - \pF \mathbf{I} \right ) \cdot \bfn}
  &= 0 & &\quad \text{on } \IFF .
  \label{eq:fluid-interface-jump-2}
\end{alignat}
Here $\jump{v} = v_1 - v_2$ denotes the jump in a function (or each
component of a vector field) $v$ over the interface $\IFF$ where $v_i
= v|_{\OF_i}$ denotes the restriction of $v$ to $\OF_i$ for $i=1,
2$. Furthermore, $\bfn$ is the unit normal of $\IFF$ directed from
$\OF_2$ into $\OF_1$ .

Correspondingly, we assume that the structure deforms as an elastic
solid satisfying the equations: find $\uS : \OS \rightarrow \R^3$ such
that
\begin{equation}
  - \nabla \cdot \sigmaS(\uS) = \bff^s \quad \text{in } \OS,
  \label{eq:structure-pde}
\end{equation}
where $\sigmaS$ is the (Cauchy) stress tensor and $\bff^s$ is a given
body force. The precise form of the Cauchy stress tensor will depend
on the choice of the elastic constitutive relation. In later sections,
we will consider both linearly elastic and hyperelastic constitutive
equations relating the displacement to the stress. As boundary
conditions, we assume that the displacement of the structure is given
on part of the boundary and that the structure experiences a boundary
traction $\bft^s_N$ on the fluid--structure interface:
\begin{alignat}{3}
  \uS &= \bfg^s_D && \quad \text{on } \partial \OS_D,
  \label{eq:structure-bc1}
  \\ \sigmaS(\uS) \cdot \bfn &= \bft^s_N && \quad \text{on } \IFS .
  \label{eq:structure-bc2}
\end{alignat}

The coupling between the fluid and the structure problems requires the
fluid and solid stresses and velocities to be in equilibrium at the
interface $\IFS$. In the stationary case considered here, these
kinematic and kinetic continuity conditions are taken care of by
ensuring that~\eqref{eq:fluid-bc1} and
\begin{equation}
  \label{eq:kinematic-matching-conditions}
  \bft^s_N = \sigmaF(\uF) \cdot \bfn
\end{equation}
hold, where $\sigmaF$ is the fluid stress tensor: $\sigmaF(\uF, \pF) =
2 \nu^f \epsilon(\uF) - \pF \mathbf{I}$ and $\epsilon(\uF)$ is the
symmetric gradient $\epsilon(\uF) = \tfrac{1}{2}(\nabla \uF + \nabla
({\uF})^{\top})$.

In summary, the stationary fluid--structure interaction problem
considered in this work is completely described by the set of
equations~\eqref{eq:fluid-momentum-pde}--\eqref{eq:kinematic-matching-conditions}.

\section{An overlapping finite element discretization of the FSI
  problem}
\label{sec:chimera-method}

The nonlinear nature of the fluid--structure interaction
problem~\eqref{eq:fluid-momentum-pde}--\eqref{eq:kinematic-matching-conditions}
mandates a nonlinear solution scheme such as a Newton-type or
fixed-point method. A classical and well-studied approach is to
decompose the coupled problem into separate systems of equations via a
Dirichlet--Neumann fixed-point
iteration~\cite{Nobile2001,LeTallec2001,Kuttler2009}. This is also the
route taken here.  Alternatively, more sophisticated iteration schemes based
on a Robin-type reformulation of the interface
conditions~\eqref{eq:fluid-bc1},\eqref{eq:structure-bc2}, and
\eqref{eq:kinematic-matching-conditions} might be employed,
see for instance \cite{Badia2008,Badia2009,BadiaNobileVergara2008}.
The basic idea of the Dirichlet--Neumann fixed-point iteration
is to start with solving the fluid
problem~\eqref{eq:fluid-momentum-pde}--\eqref{eq:fluid-interface-jump-2}
on a given starting domain. The resulting fluid boundary traction
acting on the fluid--structure interface then serves as Neumann data
for the structure
problem~\eqref{eq:structure-pde}--\eqref{eq:kinematic-matching-conditions}. The
structure deformation dictates a displacement of the fluid domain
boundary, and in turn, a new configuration of the fluid domain. This
sequence of steps is repeated until convergence.

Each of the three subproblems (the fluid problem, the structure
problem and the domain deformation) will be solved numerically using
separate finite element discretizations. Overall, we will employ an
overlapping mesh method in which a fixed background mesh is used for
part of the fluid domain and a moving mesh is used for the combination
of the structure domain and its surrounding fluid domain. We note that
methods based on overlapping meshes (as the one considered here) are
sometimes also called Chimera methods. Before describing each of the
discretizations, we here present an overview of the set-up of the
computational domains.

\begin{figure}[htb]
  \begin{center}
    \includegraphics[width=0.45\textwidth]{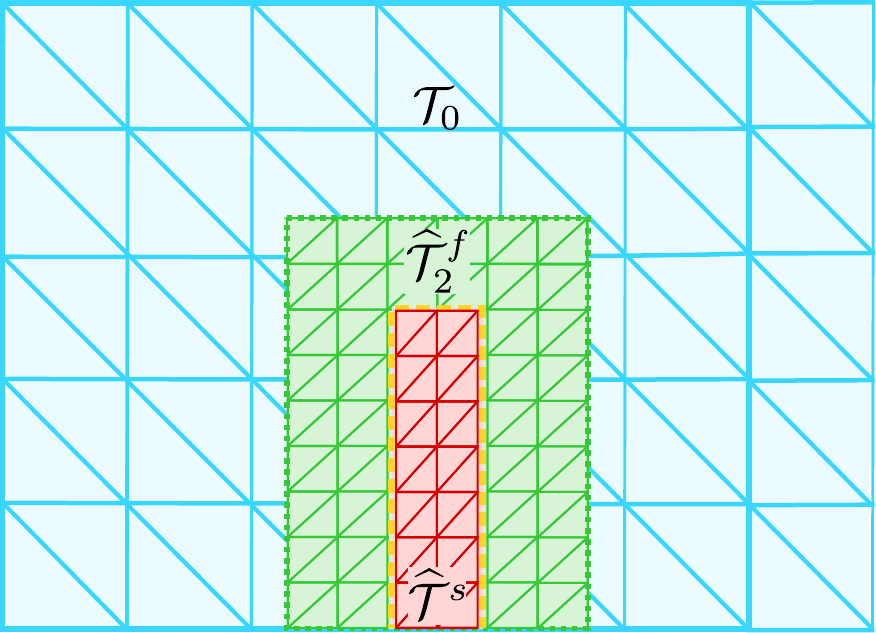}
    \hspace{1ex}
    \includegraphics[width=0.45\textwidth]{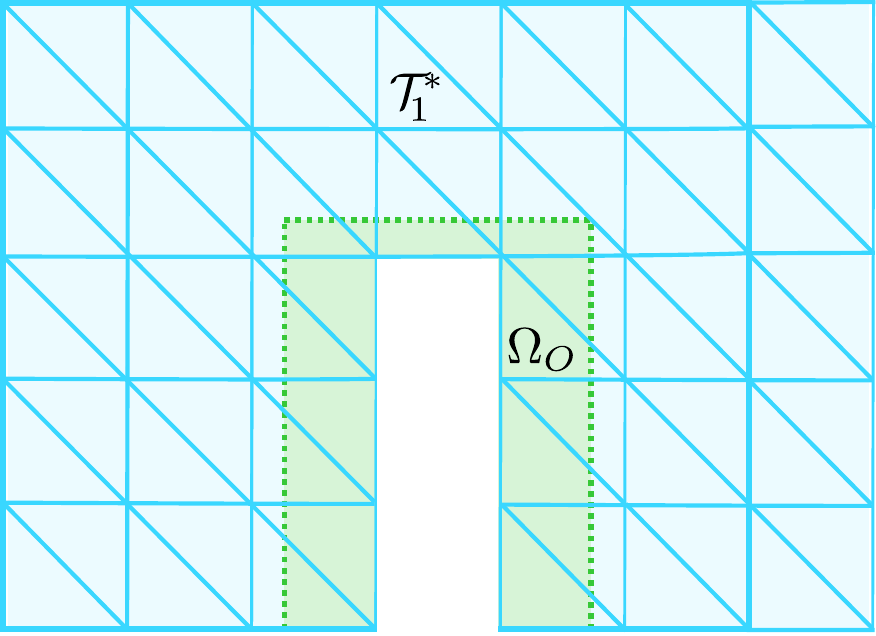}
    \caption{Chimera-mesh configuration of the computational domain in
      the starting step of the fixed-point iteration. Left: Fixed
      fluid background mesh $\mesh_0$ overlapped by the structure mesh
      $\widehat{\mesh}^s$ and a surrounding fitting fluid mesh
      $\widehat{\mesh}^f_2$. Right: Reduced fluid background mesh
      $\mesh^{\ast}_1$ and fluid overlap region
      $\Omega_O$.}  \label{fig:chimera-mesh-configuration} \end{center}
\end{figure}

For simplicity, we assume that the computational domain $\Omega$ is
fixed throughout the fixed-point iteration while the fluid and
structure subdomains will be updated in each iteration step. In each
step, we consider the following set-up, illustrated in
Figure~\ref{fig:chimera-mesh-configuration}, of the computational
domains. First, we assume that $\Omega$ is tessellated by a background
mesh $\meshO$. Second, we assume that the current representation of
the subdomains $\OF_2$ and $\OS$ are tessellated by meshes
$\fluidmesh$ and $\solidmesh$, respectively, and that these meshes
match at their common interface. As a result, $\fluidsolidmesh =
\fluidmesh \cup \solidmesh$ defines an admissible and conforming mesh
of the combined domain $\fluidsoliddomain = (\overline{\OF_2} \cup \overline{\OS})^{\circ}$. All
meshes are assumed to be admissible and to consist of shape-regular
simplices.

We further note that the background tessellation $\meshO$ may be
decomposed into three disjoint subsets:
\begin{equation}
  \label{p3:eq:mesh-splitting}
  \meshO = \mesh_{0,1} \cup \mesh_{0,2} \cup \mesh_{0, \Gamma}.
\end{equation}
Here $\mesh_{0,1}, \mesh_{0,2}, \mesh_{0,\Gamma}$ are defined with
reference to $\fluidsoliddomain$ and denote the sets of elements in
$\meshO$ that are \emph{not}, \emph{completely} or \emph{partially}
overlapped by $\fluidsoliddomain$. More precisely, $\mesh_{0,1} = \{ T
\in \mesh_0 : T \subset \overline{\OF_1} \}$, $\mesh_{0,2} = \{ T \in
\mesh_0 : T \subset \overline{\fluidsoliddomain} \}$ and
$\mesh_{0,\Gamma} = \{ T \in \mesh_0 : | T \cap \OF_1 | > 0 \text{ and
} |T \cap \fluidsoliddomain| > 0 \}$. In addition, we
assume  that $\meshO$ is sufficiently
fine near the fluid--fluid interface
in the sense that $T \cap \OS = \emptyset$ for all $T \in \mesh_{0,
  \Gamma}$.
In other words, the elements in the fluid background mesh have to be small enough
close to $\Gamma_{ff}$
such that a single element does not stretch from the fluid--fluid interface
to the fluid--structure interface.
Next, we introduce the
\emph{reduced} background mesh $\reducedmesh$, consisting of the
elements in $\meshO$ that are either not or only partially overlapped
by $\fluidsoliddomain$, and associated domain $\reduceddomain$:
\begin{align}
  \reducedmesh &= \mesh_{0,1} \cup \mesh_{0, \Gamma}, \quad
  \reduceddomain = \bigcup_{T \in \reducedmesh} T.
  \label{p3:eq:meshast-1}
\end{align}
Note that $\reduceddomain$ contains (but is generally larger than)
$\OF_1$. We further define the so-called \emph{fluid overlap region}
$\fluidoverlapregion = \OF_2 \cap \reduceddomain$. In general, for
each overlapping mesh configuration described by some (background)
mesh and some overlapping domain, the procedure described above
defines what we shall refer to as the reduced (background) mesh.

\subsection{An overlapping mesh method for the fluid problem}
\label{ssec:fluid-discretization}

In this section, we present a finite element discretization
of~\eqref{eq:fluid-momentum-pde}--\eqref{eq:fluid-interface-jump-2}
posed on a pair of overlapping meshes, first proposed
in~\citet{MassingLarsonLoggEtAl2013}. The pair of meshes consist of an overlapped
mesh and an overlapping mesh: in our case the reduced background mesh
$\reducedmesh$ plays the role of the overlapped mesh, while
$\fluidmesh$ is the overlapping mesh.

For any given mesh $\mesh$, let $V_h(\mesh)$ be the space of
continuous piecewise linear vector fields and let $Q_h(\mesh)$ be the
space of continuous piecewise linears, both defined relative to
$\mesh$. We define the composite finite element spaces $V_h$ and $Q_h$
for the overlapping fluid meshes by
\begin{equation}
  \label{p3:eq:elementspaces}
  V^f_h = V_h(\reducedmesh) \bigoplus V_h(\fluidmesh), \quad
  Q^f_h = Q_h(\reducedmesh) \bigoplus Q_h(\fluidmesh).
\end{equation}
Moreover, we denote by $V^f_{h, \bfg^f}$ the subspace of $V^f_h$ that
satisfies the boundary
conditions~\eqref{eq:fluid-bc1}--\eqref{eq:fluid-bc2} and by $V^f_{h,
  \boldsymbol{0}}$ the corresponding homogeneous version.  The
overlapping mesh discretization
of~\eqref{eq:fluid-momentum-pde}--\eqref{eq:fluid-interface-jump-2}
then reads: find $(\uF_h, \pF_h) \in V^f_{h, \bfg^f} \times Q^f_h$
such that
\begin{equation}
  \label{eq:stokes-olm}
  A^f_h(\uF_h, \pF_h ; \bfv, q) = L^f_h(\bfv, q)
  \quad \foralls (\bfv, q) \in V^f_{h, \boldsymbol{0}} \times Q^f_h,
\end{equation}
where $A^f_h$ is defined for all $\bfu, \bfv \in V^f_h$ and all $p, q
\in Q^f_h$ by
\begin{equation}
  \begin{split}
    \label{eq:Ah-olm}
    A^f_h(\bfu, p; \bfv, q) = a^f_h(\bfu, \bfv) + b^f_h(\bfv, p) +
    b^f_h(\bfu, q) + i^f_h(\bfu, \bfv) - j^f_h(p, q),
  \end{split}
\end{equation}
and the forms $a^f_h$, $b^f_h$, $i^f_h$ and $j^f_h$ are given by
\begin{align}
  \label{eq:ah-olm}
  a^f_h(\bfu, \bfv) &= (\nabla \bfu , \nabla{\bfv})_{\OF_1 \cup \OF_2}
  - (\meanvalue{\nablan \bfu}, \jump{{\bfv}} )_{\IFF} -
  (\meanvalue{\nablan \bfv}, \jump{{\bfu}} )_{\IFF}
  + \gamma (h^{-1} \jump{{\bfu}}, \jump{{\bfv}} )_{\IFF}, \\
  \label{eq:bh-olm}
  b^f_h(\bfv, q) &= - (\nabla \cdot \bfv, q)_{\OF_1 \cup \OF_2}
  + (\jump{\bfv} \cdot \bfn, \meanvalue{q})_{\IFF}, \\
  \label{eq:ih-olm}
  i^f_h(\bfu, \bfv)
  &= (\nabla(\bfu_1 - \bfu_2), \nabla(\bfv_1 - \bfv_2))_{\fluidoverlapregion}, \\
  \label{eq:jh-olm}
  j^f_h(p, q) &= \delta \sum_{T \in \reducedmesh \cup \fluidmesh }
  h_T^2 (\nabla p, \nabla q)_T,
\end{align}
for $\delta > 0$. Here and throughout, $(\cdot, \cdot)_{K}$ denotes
the $L^2(K)$ inner product over some domain $K$, while $\meanvalue{v}$
denotes a convex combination $\meanvalue{v} = \alpha_1 v_1 + \alpha
v_2$ with $\alpha_1 + \alpha_2 = 1$ of $v$ across the interface
$\IFF$. In particular, we choose $\meanvalue{v} = v_2$ in accordance
with~\citet{HansboHansboLarson2003}. Finally, the linear form $L^f_h$
is defined by
\begin{equation}
  L^f_h(\bfv, q) = (\bff^f, \bfv)
  - \delta \sum_{T\in \reducedmesh \cup \fluidmesh}
  h_T^2 (\bff^f, \nabla q)_T
  \label{eq:Lh-olm}
\end{equation}
for all $\bfv \in V^f_h$ and all $q \in Q^f_h$.

A major strength of the employed scheme for the fluid problem is that
the extension of the stabilization term~\eqref{eq:jh-olm} from the
physical domain $\Omega^f_1$ to the overlap region
$\fluidoverlapregion$ in combination with the least-square
stabilization~\eqref{eq:ih-olm} results in a well-conditioned and
optimally convergent scheme, independent of the location of the
overlapping mesh with respect to the fixed background mesh. Thereby,
typical difficulties arising from potentially small cut cells where
$|T \cap \Omega^f_2 | \ll | T |$ for $T \in \mathcal{T}_{0,\Gamma}$
are completely eliminated. Consequently, for a continuous solution
$(\uF, \pF)$ satisfying
of~\eqref{eq:fluid-momentum-pde}--\eqref{eq:fluid-interface-jump-2}
and a discrete solution $(\uF_h, \pF_h)$
satisfying~\eqref{eq:stokes-olm}, the following optimal error estimate
holds independently of the fluid--fluid interface
position~\citep{MassingLarsonLoggEtAl2013}:
\begin{equation}
  \tn ( \uF - \uF_h , \pF - \pF_h ) \tn
  \leqslant C h \left( | \uF |_{2, \OF} + |\pF|_{1, \OF} \right).
  \label{eq:a-priori}
\end{equation}
Here, $\tn \cdot \tn$ is an appropriate version of the standard norm
on $H^1(\OF) \times L^2(\OF)$ accounting for the fluid overlap region
$\fluidoverlapregion$; see~\citep{MassingLarsonLoggEtAl2013} for more details.

\subsection{A finite element discretization of the structure problem}

The structure problem is described
by~\eqref{eq:structure-pde}--\eqref{eq:structure-bc2} in the current
solid domain. As the current solid domain is actually unknown, a
standard approach to discretizing such problems is to map the
governing equations back to a fixed reference (Lagrangian) frame. We
choose a reference domain $\refe{\Omega}^s$ with coordinates
$\refe{\bfx}$ and denote the deformation map from the reference to the
current solid domain by $\bphi^s$:
\begin{equation}
  \bfx = \bphi^s(\refe{\bfx})
  \quad \text{for } \refe{\bfx} \in \refe{\Omega}^s .
\end{equation}
In general, the notation for all domains and quantities pulled back to
the Lagrangian framework will be endowed with a $\refe{\phantom{u}}$;
for instance $\refe{\Omega}^s$ and $\refe{\bfu}^s$ denote the solid
reference domain and solid displacement in the reference frame,
respectively. In particular, $\boldsymbol \phi^s = \mathbf{I} +
\refe{\bfu}^s$.

In the Lagrangian frame, the problem reads: find the solid
displacement $\refe{\bfu}^s: \refe{\Omega}^s \to \R^3$ such that
\begin{alignat}{3}
  \label{eq:structure-pde-ref}
  - \nabla \cdot \refe{\Pi} (\refe{\bfu}^s) &= \refe{\bff}^s
  && \quad \text{in } \refe{\Omega}^s, \\
  \label{eq:structure-bc1-ref}
  \refe{\bfu}^s &= \refe{\bfg}^s_D
  && \quad \text{on } \partial \refe{\Omega}^s_D, \\
  \label{eq:structure-bc2-ref}
  \refe{\Pi} (\refe{\bfu}^s) \cdot \refe{\bfn} &= \refe{\bft}^s_N &&
  \quad \text{on } \refe{\Gamma}^{fs} .
\end{alignat}
Here, the displacement $\refe{\bfu}^s$ and the boundary displacement
$\refe{\bfg}^s_D$ result from the standard affine pull-back of the
corresponding quantities in the current domain, for instance
$\refe{\bfu}^s(\refe{\bfx}) = \bfu^s(\bfx)$, and $\refe{\bfn}$ is the
outward normal of the fluid--structure interface in the reference
frame. Further, let $\bfF^s = \nabla \bphi^s$ and $J^s = \det
\bfF^s$. We let $\refe{\bff}^s(\refe{\bfx}) = J^s
\bff^s(\bfx)$. Moreover, $\refe{\Pi}(\refe{\bfu}^s)$ denotes the first
Piola-Kirchhoff stress tensor, resulting from a Piola transformation
of the Cauchy stress tensor $\sigmaS$:
\begin{equation}
  \refe{\Pi}(\refe{\bfu}^s) (\refe{\bfx})
  = J^s(\refe{\bfx}) \,
  \sigmaS (\bphi^s( \refe{\bfx})) (\bfF^s)^{-\top} (\refe{\bfx}) .
\end{equation}
In view of~\eqref{eq:kinematic-matching-conditions}, we will enforce
that the boundary traction acting on the solid in the reference domain
is the Piola transform of the fluid traction exerted on the
fluid--structure interface by the fluid in the current or physical
configuration. This will be detailed in Section~\ref{sec:algorithm}.

The governing
equations~\eqref{eq:structure-pde-ref}--\eqref{eq:structure-bc2-ref}
must be completed by a constitutive equation relating the stress to
the strain. In the case of a hyperelastic material, by definition,
there exists a strain energy density $\Psi$ such that
\begin{equation}
  \refe{\Pi}(\bfF) = \frac{\partial \Psi}{\partial \bfF} .
\end{equation}
One example is the St.~Venant--Kirchhoff material model, in which
\begin{equation}
  \Psi(\bfF) = \mu \tr \bfE^2 + \frac{\lambda}{2} (\tr \bfE)^2,
  \quad
  \text{where } \bfE = \frac{1}{2} (\bfF^{\top} \bfF - \mathbf{I}),
  \label{eq:venant-kirchhoff-material}
\end{equation}
for Lam\'e constants $\mu, \lambda > 0$.

In the special case of a linearly elastic material, we assume that the
reference and physical configurations coincide, so
that~\eqref{eq:structure-pde}--\eqref{eq:structure-bc2} hold over
$\refe{\Omega}^s$ directly with $\sigmaS(\uS) = 2 \mu \epsilon(\uS) +
\lambda \tr (\epsilon(\uS)) \mathbf{I}$.

To solve~\eqref{eq:structure-pde-ref}--\eqref{eq:structure-bc2-ref}
numerically, let $\refe{\mesh}^s$ be a tessellation of
$\refe{\Omega}^s$ such that $\solidmesh = \bphi^s(\refe{\mesh}^s)$ and
introduce the finite element approximation space
\begin{equation}
  \refe{V}^s_{h, \bfg}
  = \{ \bfv \in V_h(\refe{\mesh}^s) \quad
  \text{such that} \quad \bfv|_{\partial \refe{\Omega}^s_D} = \bfg \},
\end{equation}
where $V_h(\refe{\mesh}^s)$ is the space of continuous piecewise
linear vector fields defined relative to $\refe{\mesh}^s$ as
before. The finite element formulation
of~\eqref{eq:structure-pde-ref}--\eqref{eq:structure-bc2-ref} then
reads: find $\refe{\bfu}^s_h \in \refe{V}^s_{h, \refe{\bfg}^s_D}$ such
that
\begin{align}
  (\refe{\Pi}(\refe{\bfu}^s_h), \nabla \bfv)_{\refe{\Omega}^s} -
  (\refe{\bft}^s_N, \bfv)_{\refe{\Gamma}^{fs}} - (\refe{\bff}^s,
  \bfv)_{\refe{\Omega}^s} &= 0 \quad \foralls \bfv \in \refe{V}^s_{h,
    \boldsymbol{0}}.
  \label{eq:structure-problem-weak-form}
\end{align}
Note that the generally nonlinear constitutive relation and the
geometric nonlinearity mandate a nonlinear solution scheme, such as a
Newton method or an inner fixed-point iteration
for~\eqref{eq:structure-problem-weak-form}.

\subsection{Deformation of the surrounding fluid domain}
\label{ssec:mesh-motion}

The overlapping mesh method relies on keeping the background part of
the fluid domain $\OF_1$ fixed while moving the part of the fluid
domain $\OF_2$ surrounding the structure. This movement ensures that
the mesh $\fluidmesh$ of the latter part of the fluid domain and the
structure mesh $\solidmesh$ match at the fluid--structure
interface. The movement is dictated by the structure deformation only
at the fluid--structure interface: the motion of the interior of the
fluid domain $\OF_2$ is subject to numerical modeling. Standard
approaches for the domain motion include mesh smoothing via
diffusion-type equations or treating the fluid domain as a
pseudo-elastic structure. Here, we choose the latter approach and
model the deformation of the fluid domain as a linearly elastic
structure. This approach allows for typically larger deformations than
a simple diffusion equation based mesh smoothing, while avoiding
unnecessary complexity.

We start with a fixed reference domain $\refe{\Omega}^f_2$ and
consider the following mesh deformation problem over this domain: find
the mesh displacement $\refe{\bfu}^m : \refe{\Omega}^f_2 \rightarrow
\R^3$ such that
\begin{alignat}{3}
  - \nabla \cdot \refe{\bfsigma}^m (\refe{\bfu}^m) &= 0 \quad &&
  \text{in } \refe{\Omega}^f_2,
  \label{eq:mesh-motion-strong}
  \\
  \refe{\bfsigma}^m (\refe{\bfu}^m) \cdot \refe{\bfn} &= 0 \quad
  &&\text{on } \refe{\Gamma}^{ff},
  \label{eq:mesh-motion-nbc}
  \\
  \refe{\bfu}^m &= \refe{\bfu}^s \quad &&\text{on }
  \refe{\Gamma}^{fs},
  \label{eq:mesh-motion-dbc}
\end{alignat}
where the stress tensor $\refe{\bfsigma}^m$ is given by
\begin{equation}
  \refe{\bfsigma}^m(\refe{\bfu}^m)
  = 2 \mu_m \epsilon(\refe{\bfu}^m) +
  \lambda_m \tr(\epsilon (\refe{\bfu}^m)) \mathbf{I}
  \label{eq:hookes-law}
\end{equation}
for chosen Lam\'e constants $\mu_m, \lambda_m > 0$.
Let now
$\refe{\mesh}^f_2$ be a tessellation of $\refe{\Omega}^f_2$. We define
the finite element space $\refe{V}^m_{h, \bfg}$ by
\begin{equation}
  \refe{V}^m_{h, \bfg}
  = \{ \bfv \in V_h(\refe{\mesh}^f_2) \quad
  \text{such that} \quad \bfv|_{\refe{\Gamma}^{fs}} = \bfg \}.
\end{equation}
The corresponding finite element formulation of the mesh
problem~\eqref{eq:mesh-motion-strong}--\eqref{eq:mesh-motion-dbc} is
then: find $\refe{\bfu}^m_h \in \refe{V}^m_{h, \refe{\bfu}^s}$ such
that
\begin{equation}
  (\refe{\bfsigma}^m (\refe{\bfu}^m_h), \bfv)_{\refe{\Omega}^f_2} =
  0 \quad \foralls \bfv \in \refe{V}^m_{h, \boldsymbol{0}}.
\end{equation}

Finally, we define $\mesh^{f}_2 = \bphi^m_h(\refe{\mesh}^f_2)$ with
the discrete mesh deformation $\bphi^m_h = \mathbf{I} +
\refe{\bfu}^m_h$.
The current surrounding fluid domain is then defined
accordingly: $\Omega^{f}_2 = \bphi^m_h(\refe{\Omega}^f_2)$.
The use of boundary condition~\eqref{eq:mesh-motion-nbc}
ensures that the fluid--structure interface is
preserved in the sense that
\begin{align}
\Gamma^{fs} = \partial \Omega^{f}_2 \cap
\partial \Omega^{s}
= \bphi^m_h(\refe{\Gamma}^{fs})
= \bphi^s_h(\refe{\Gamma}^{fs}).
\end{align}
where $\bphi_h^s$ is the solid deformation
given by the discrete solution $\refe{\bfu}^{s}_h$
of problem~\eqref{eq:structure-problem-weak-form}.

\section{Solution algorithm for the discretized FSI problem}
\label{sec:algorithm}

We are now in a position to give a detailed description of the overall
solution scheme for the fully coupled fluid--structure interaction
problem. We start with reviewing the formulation of the
fluid--structure coupling in the discrete setting. For the discrete
formulation, a third interface condition~\eqref{eq:mesh-motion-dbc}
needs to be added to the two interface conditions~\eqref{eq:fluid-bc1}
and ~\eqref{eq:structure-bc2}, due to the additional mesh deformation
problem described in Section~\ref{ssec:mesh-motion}.  The mesh
deformation allows to express the fluid stress tensor acting on $\IFS$
in the reference configuration $\refe{\Gamma}^{fs}$ via a Piola
transformation. Consequently, the stress equilibrium
condition~\eqref{eq:structure-bc2} at the fluid--structure interface
can be reformulated in the Lagrangian frame according
to~\eqref{eq:structure-bc2-ref}. In summary, the discrete formulation
of the fluid--structure interface conditions reads:
\begin{align}
  \label{eq:final-fsi-condition-1}
  \uF &= 0 \quad \;\;\;\, \text{on } \IFS,
  \\
  \label{eq:final-fsi-condition-2}
  \refe{\bfu}^s &= \refe{\bfu}^m \quad \text{on } \refe{\Gamma}^{fs},
  \\
  \label{eq:final-fsi-condition-3}
  \refe{\Pi}(\refe{\bfu}^s) (\refe{\bfx}) \cdot
  \refe{\bfn}(\refe{\bfx}) &= J^m(\refe{\bfx}) \, \sigmaF (\bphi^m(
  \refe{\bfx})) (\bfF^m)^{-\top} (\refe{\bfx}) \cdot
  \refe{\bfn}(\refe{\bfx}) \quad \text{on } \refe{\Gamma}^{fs}.
\end{align}
As outlined in Section~\ref{sec:chimera-method}, we employ a classical
Dirichlet--Neumann fixed-point iteration approach to ensure that the
interface
conditions~\eqref{eq:final-fsi-condition-1}--\eqref{eq:final-fsi-condition-3}
are approximately satisfied by the computed solution within a user
provided tolerance.  The iteration scheme is presented in detail in
Algorithm~\ref{alg:fixed-point-iteration}, where the relaxation
parameter $\omega^{i}$ was chosen dynamically to accelerate the
convergence of the fixed-point iteration.  Moreover, the fluid
boundary traction is incorporated as Neumann data in the weak
formulation of the structure problem by a properly chosen functional
representing the boundary traction weighted with some given test
function.  A thorough explanation of both of these intermediate steps
will be given in the next sections.

\begin{algorithm}[htb]
  \begin{tabbing}
    \\
    $\refe{\bfu}^{s, k} := 0$
    \\
    $\refe{\bfu}^{m, k} := 0$
    \\
    \textbf{do }
    \\
    \tab \textbf{Update overlapping fluid meshes}
    \\
    \tab $\Omega^{s,k+1} := (\mathbf{I} + \refe{\bfu}^{s,
      k})(\refe{\Omega}^s)$
    \\
    \tab $\Omega^{f,k+1}_2 := (\mathbf{I} + \refe{\bfu}^{m,
      k})(\refe{\Omega}^f_2)$
    \\
    \tab $\Omega^{fs,k+1} :=  \Omega^{s,k+1} \cup \Omega^{f,k+1}_2$
    \\
    \tab Compute reduced background mesh $(\mesh^{f,k+1}_1)^{\ast}$
    with respect to $\Omega^{fs,k+1}$
    \\
    \tab $\meshFk := (\mesh^{f,k+1}_1)^{\ast} \cup \meshFtwok$
    \\
    \\
    \tab \textbf{Solve fluid problem}
    \\
    \tab Find $(\bfu^{f, k+1}_h,p^{f, k+1}_h)$ such that $\foralls (\bfv^{f, k+1}_h,q^{f, k+1}_h) \in
    V_h^{f,k+1} \times Q_h^{f,k+1}$
  \end{tabbing}
  \[
  A^{f,k}_h(\bfu^{f, k+1}_h,p^{f, k+1}_h;\bfv^{f, k+1}_h, q^{f,
    k+1}_h) = L^{f,k+1}(\bfv^{f, k+1}_h,q^{f, k+1}_h)
  \]
  \begin{tabbing}
    \\
    \tab \textbf{Update boundary traction functional}
    \\
    \tab Define $L^{fs,k+1}(\cdot)$ by
  \end{tabbing}
  \[
  L^{fs,k+1}(\refe{\bfv}^{s, k+1}_h) := R^{f,k+1}(\bfu^{f,
    k+1}_h,p^{f, k+1}_h;\bfv^{f, k+1}_h)
  \]
  \begin{tabbing}
    \\
    \tab \textbf{Solve structure problem}
    \\
    \tab Find $\refe{\bfu}^{s, k+1}_h$ such that $\foralls \refe{\bfv} \in \widehat{V}_h^{s}$
  \end{tabbing}
  \[
  A^{s}_h(\refe{\bfu}^{s, k+1}_h, \refe{\bfv}) = L^{s}(\refe{\bfv}) +
  L^{fs,k+1}(\refe{\bfv})
  \]
  \begin{tabbing}
    \\
    \tab \textbf{Dynamic relaxation}
    \\
    \tab Compute  $\omega^{k+1}$ according to~\eqref{eq:compute-relax-parameter}
    \\
    \tab $\refe{\bfu}^{s, k+1}_h
    := \omega^{k+1} \refe{\bfu}^{s, k+1}_h + (1-\omega^{k+1}) \refe{\bfu}^{s, k}_h$
    \\
    \\
    \tab \textbf{Solve mesh problem}
    \\
    \tab Find $\refe{\bfu}^{m, k+1}_h$ such that $\foralls \refe{\bfv} \in \widehat{V}_h^{m}$
  \end{tabbing}
  \begin{align*}
    A^{m}_h(\refe{\bfu}^{m, k+1}_h, \refe{\bfv}) &= L^{s}(\refe{\bfv})
    \\
    \refe{\bfu}^{m, k+1}_h &= \refe{\bfu}^{s, k+1}_h \text{on }
    \refe{\IFS}
  \end{align*}
  \begin{tabbing}
    \textbf{while } $\| \refe{\bfu}^{s, k+1}_h - \refe{\bfu}^{s, k}_h \| \leqslant $ \texttt{TOL}
  \end{tabbing}
  \caption{Fixed-point iteration.}
  \label{alg:fixed-point-iteration}
\end{algorithm}

\subsection{Dynamic relaxation}

Let $U_{_S}^k$ denote the coefficient vector of the finite element
approximation $\refe{\bfu}^{s,k}_h$ of
\eqref{eq:structure-problem-weak-form} computed in the $k$-th
iteration step.  To accelerate the convergence of the iteration
scheme, a relaxation step is introduced:
\begin{align}
  U_{_S}^{k+1} := \omega_k U_{_S}^{k+1} + (1-\omega_k) U_{_S}^{k},
  \label{eq:def:relaxation}
\end{align}
where the relaxation parameter $\omega_k$ is dynamically chosen in
each iteration step. Here, we employed Aiken's method
\cite{Kuttler2008,Kuttler2009} which is a simple scheme, yet it can
greatly improve the convergence rate compared to a fix choice of
$\omega_k$, as demonstrated by \citet{Kuttler2008,Kuttler2009}.
Introducing the residual displacement $\Delta^k U_{_S}$ by
\begin{align}
  \Delta^k U_{_S} := U_{_S}^k - U_{_S}^{k-1},
  \label{eq:def-displacement-residual}
\end{align}
the new relaxation parameter $\omega_{k+1}$ is then computed by
\begin{align}
  \omega_{k} = \max \left\{\omega_{\max}, \omega_{k-1} \left( 1 -
      \dfrac{\Delta^{k+1} U_{_S}}{\| \Delta^{k+1} U_{_S} - \Delta^{k}
        U_{_S}\|^2} \right) \right\},
  \label{eq:compute-relax-parameter}
\end{align}
where $\omega_{\max}$ is a safety parameter chosen to avoid too large
over-relaxation. The convergence of the fixed-point iteration might be
accelerated further by employing more sophisticated schemes based on
Robin--Robin coupling~\cite{Badia2008,Badia2009} or vector
extrapolation~\cite{Kuttler2009}.

\subsection{Computation of the boundary traction}

Given the solution $\uF$ and a pressure solution $p^f$ of the fluid
subproblem~\eqref{eq:fluid-momentum-pde}--\eqref{eq:fluid-bc2}, the
incorporation of the fluid boundary traction into the weak formulation
of the structure problem~\eqref{eq:structure-problem-weak-form}
requires the evaluation of the so-called weighted fluid boundary
traction on $\IFS$ defined by
\begin{align}
  L^{fs}(\bfv) &= (\sigmaF(\uF,p^f)\cdot \bfn, \bfv)_{\IFS},
  \label{eq:fluid-traction-strong}
\end{align}
for test functions $\bfv \in V^{s}$. The
functional~\eqref{eq:fluid-traction-strong} possesses various
equivalent representations in the continuous case which are no longer
equivalent when fluid velocity $\uF$ and pressure $\pF$ and test
function $\bfv$ are replaced by their discrete counterparts $\uF_h$,
$\pF_h$ and $\bfv_h \in V^s_h(\Omega)$, respectively. It has been
observed by \citet{Dorok1995,John1997,Giles1997} that
using~\eqref{eq:fluid-traction-strong} directly might lead to an
inaccurate evaluation of the weighted boundary traction. In our work,
we therefore employ an alternative formulation of the weighted
boundary traction in the form
\begin{align}
  L^{fs}(\bfv_h) &= (\sigmaF(\uF_h,\pF_h),\Ext(\bfv_h))_{\OF} -
  (\bff^f, \Ext(\bfv_h))_{\OF},
  \label{eq:fluid-traction-functional}
\end{align}
which was proposed and investigated by \citet{Giles1997} in the
context of a posteriori error estimation. Here, $\Ext(\bfv)$ is any
function in $H^1(\Omega^{fs})$ such that $\Ext(\bfv_h)|_{\IFS} =
\bfv_h$.  Compared to the naive evaluation
using~\eqref{eq:fluid-traction-strong}, the
formulation~\eqref{eq:fluid-traction-functional} was shown to compute
the weighted boundary traction more accurately and to greatly improve
the convergence of stress related quantities such as the lift and drag
coefficients.

\section{Numerical results}

We conclude this paper with two numerical tests, both in three
spatial dimensions. The numerical experiments were carried out using
the \rm{DOLFIN-OLM} library (\url{http://launchpad.net/dolfin-olm}).
We first study the convergence rates for the finite element
approximations of the fluid velocity, fluid pressure and structure
displacement by constructing an artificial fluid--structure
interaction problem possessing an analytical solution. Second, we
consider the flow around an elastic flap immersed in a
three-dimensional channel.

\subsection{Software for overlapping mesh variational formulations}

The assembly of finite element tensors corresponding to standard
variational formulations on conforming, simplicial meshes, such
as~\eqref{eq:structure-problem-weak-form}, involves integration over
elements and possibly, interior and exterior facets. In contrast, the
assembly of variational forms defined over overlapping meshes, such
as~\eqref{eq:ah-olm}-\eqref{eq:jh-olm} and~\eqref{eq:Lh-olm},
additionally requires integration over cut elements and cut
facets. These mesh entities are of polyhedral, but otherwise
arbitrary, shape. As a result, the assembly process is highly
non-trivial in practice and requires additional geometry related
preprocessing, which is challenging in particular for
three-dimensional meshes.

As part of this work, the technology required for the automated
assembly of general variational forms defined over overlapping meshes
has been implemented as part of the software library
\rm{DOLFIN-OLM}. This library builds on the core components of the
FEniCS Project~\citep{LoggMardalEtAl2011,Logg2007}, in particular
DOLFIN~\citep{LoggWells2010a}, and the computational geometry
libraries \rm{CGAL}~\citep{cgal} and
\rm{GTS}~\citep{gts}. \rm{DOLFIN-OLM} is open source and freely
available from \url{http://launchpad.net/dolfin-olm}.

There are two main challenges involved in the implementation: the
computational geometry and the integration of finite element
variational forms on cut cells and facets. The former involves
establishing a sufficient topological and geometric description of the
overlapping meshes for the subsequent assembly process. To this end,
\rm{DOLFIN-OLM} provides functionality for finding and computing the
intersections of triangulated surfaces with arbitrary simplicial
background meshes in three spatial dimensions; this functionality
relies on the computational geometry libraries \rm{CGAL} and
\rm{GTS}. These features generate topological and geometric
descriptions of the cut elements and facets. Based on this
information, quadrature rules for the integration of fields defined
over these geometrical entities are produced. The computational
geometrical aspect of this work extends, but shares many of the
features of, the previous work~\citep{Massing2012a}, and is described
in more detail in the aforementioned reference.

Further, by extending some of the core components of the FEniCS
Project, in particular FFC~\citep{KirbyLogg2006,LoggOelgaardEtAl2011a}
and UFC~\citep{AlnaesLoggEtAl2012a}, this work also provides a finite
element form compiler for variational forms defined over overlapping
meshes. Given a high-level description of the variational formulation,
low-level C++ code can be automatically generated for the evaluation
of the cut element, cut facet and surface integrals, in addition to
the evaluation of integrals over the standard (non-cut) mesh
entities. The generated code takes as input appropriate quadrature
points and weights for each cut element or facet; these are precisely
those provided by the \rm{DOLFIN-OLM} library.

As a result, one may specify variational forms defined over finite
element spaces on overlapping meshes in high-level \rm{UFL} notation
\citep{AlnaesEtAl2012,Alnaes2011a}, define the overlapping fluid meshes $\{\mesh_0,
\mesh^f_2\}$ and then invoke the functionality provided by the
\rm{DOLFIN-OLM} library to automatically assemble the corresponding
stiffness matrix. In particular, the numerical experiments presented
below, employing the variational formulation defined
by~\eqref{eq:stokes-olm}, have been carried out using this technology.

\subsection{Convergence test}
While numerical studies presented in \cite{MassingLarsonLoggEtAl2013} confirmed the
theoretically predicted convergence rates for the overlapping mesh
method for the pure flow problem presented in
Section~\ref{ssec:fluid-discretization}, we here conduct a convergence
study of the coupled FSI problem to verify the overall solution
algorithm as described in Algorithm~\ref{alg:fixed-point-iteration}.
To examine the convergence rates for the finite element approximations
of the fluid velocity, fluid pressure and structure displacement, we
construct a stationary FSI problem with a known analytical solution by
employing the method of manufactured solutions as outlined in the
following. The detailed analytical derivation of the fluid and
structure related quantities are not included here to keep the
presentation at an appropriate length, but can be obtained as an
IPython based notebook available at
\url{http://nbviewer.ipython.org/6291921}.

\begin{figure}[htb]
  \centering
  \includegraphics[width=0.45\textwidth]{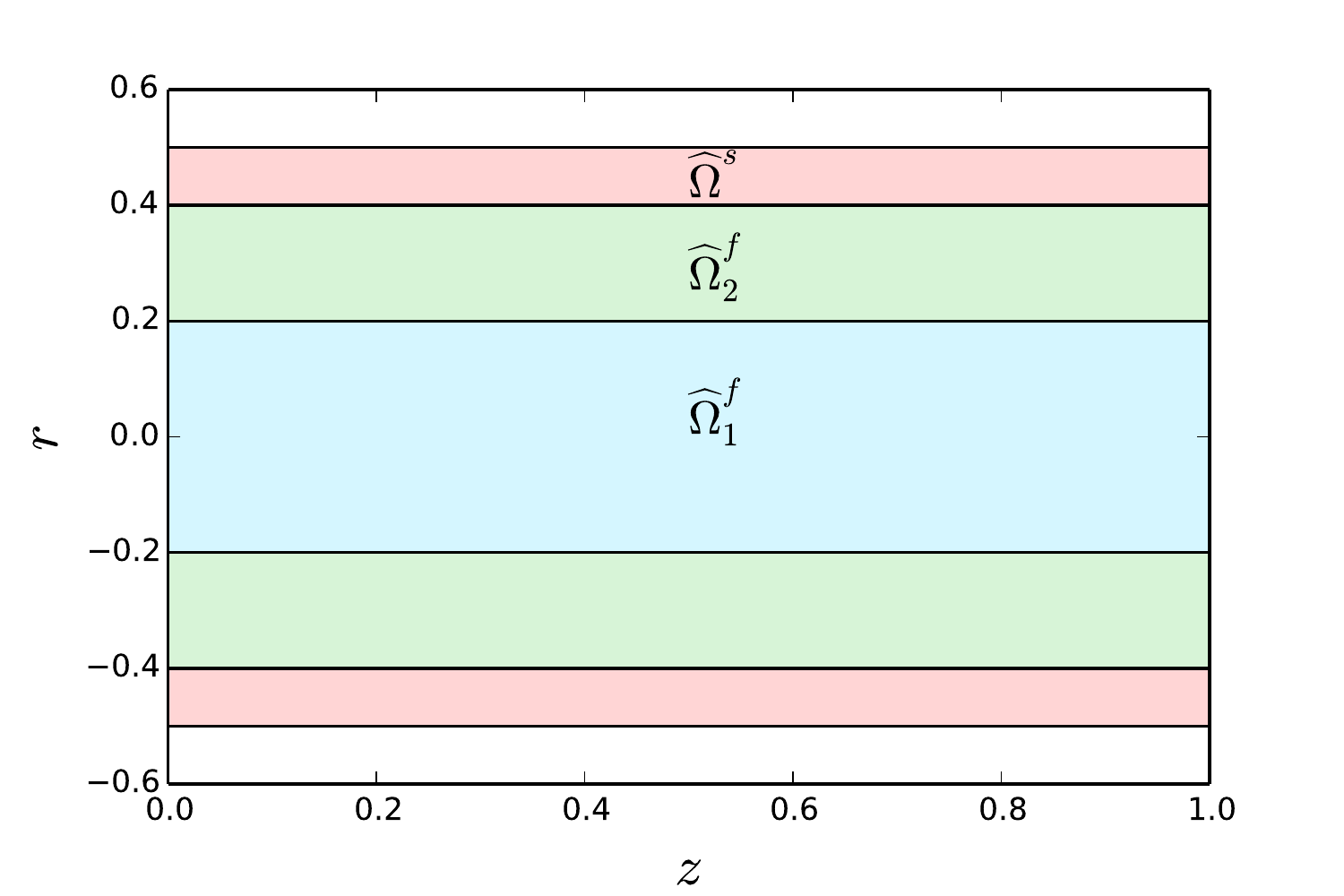}
  \includegraphics[width=0.45\textwidth]{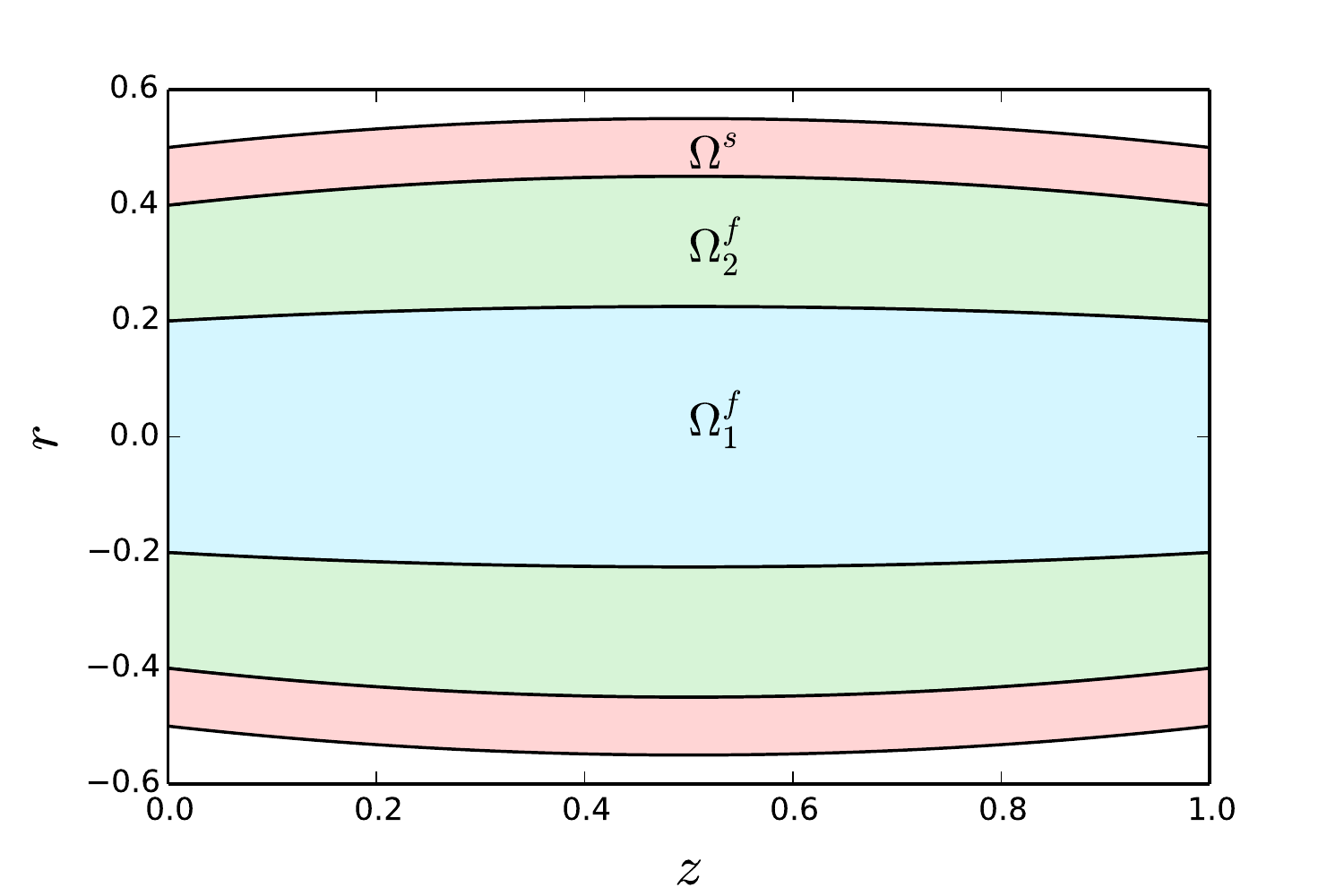}
  \caption{Cross-section through the cylinder-symmetric reference
    (left) and physical (right) domains for the analytical FSI
    reference problem.}
  \label{fig:analytical-fsi-example-domains}
\end{figure}

In the reference configuration, the fluid domain $\refe{\Omega}^f$
consists of a straight tube of length $L = 1.0$ and diameter $R^f =
0.4$.  We decompose $\refe{\Omega}^f$ into into a tube of radius
$R^f_1 = 0.3$ and a cylinder annulus satisfying $0.3 \leqslant r
\leqslant 0.4 = R^f_2$. The solid domain $\refe{\Omega}^s$ is given by
a cylinder annulus of thickness $H^s = 0.1$ surrounding the fluid
domain $\refe{\Omega}^f$. Using cylinder coordinates, the displacement
$\refe{\bfu}^s$ of the solid domain is prescribed by a purely radial,
$z$-dependent translation
\begin{equation}
  \refe{\bfu}^s(r,\varphi,z) = H(z)\boldsymbol{e}_r,
  \label{eq:solid-radial-translation}
\end{equation} where $H(z) = H^s 2 z (1-z)$.
Correspondingly, the deformation of the fluid domain is determined by
a radial stretching of the form
\begin{equation}
  \refe{\bfu}^m(r,\varphi,z) = \rho(1 + H(z)/R^f)\boldsymbol{e}_r.
  \label{eq:fluid-radial-stretching}
\end{equation}
The reference and physical configuration of the various domains are
depicted in Figure~\ref{fig:analytical-fsi-example-domains}.

To obtain a divergence-free velocity field in the final physical
configuration, the fluid velocity is defined as a simple parabolic
channel flow on the reference domain and then mapped to the physical
domain via the Piola transformation induced by the fluid domain
deformation~\eqref{eq:fluid-radial-stretching}.  For the pressure, we
simply choose $p(x,y,z) = 1-z$. Since the interface
condition~\eqref{eq:final-fsi-condition-3} is not satisfied exactly,
we introduce an auxiliary traction $\bft_a$ given by the non-vanishing
jump in the normal stresses:
\begin{equation}
  \bft_a
  = \bigl( \refe{\Pi}(\refe{\bfu}^s) (\refe{\bfx})
  - J^m(\refe{\bfx}) \,
  \sigmaF (\bphi^m( \refe{\bfx})) (\bfF^m)^{-\top} (\refe{\bfx}) \bigr)\cdot
  \refe{\bfn}^s
  \quad \text{on } \refe{\Gamma}^{fs}.
  \label{eq:final-fsi-condition-3-with-auxiliary-traction}
\end{equation}
Regarding the remaining boundary parts, the solid displacement is
uniquely determined by imposing the given displacement $\refe{\bfu}^s$
as a Dirichlet boundary condition on $\partial\refe{\Omega}^s
\setminus \refe{\Gamma}^{fs}$. For the fluid problem, we prescribed
the velocity profile on the inlet and impose the zero pressure on the
outlet.

\begin{figure}[htb]
  \centering
  \includegraphics[width=0.45\textwidth]{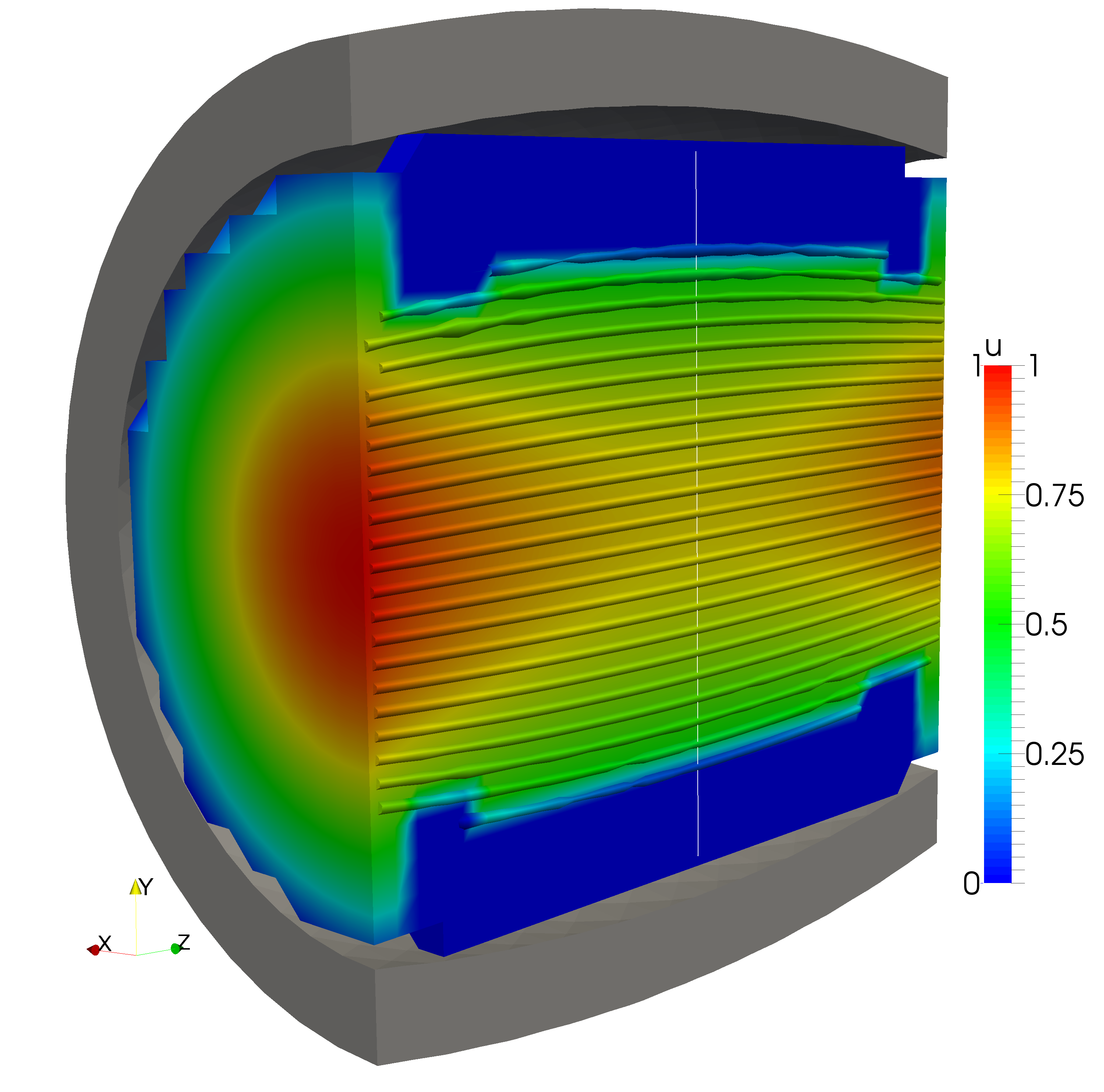}
  \includegraphics[width=0.45\textwidth]{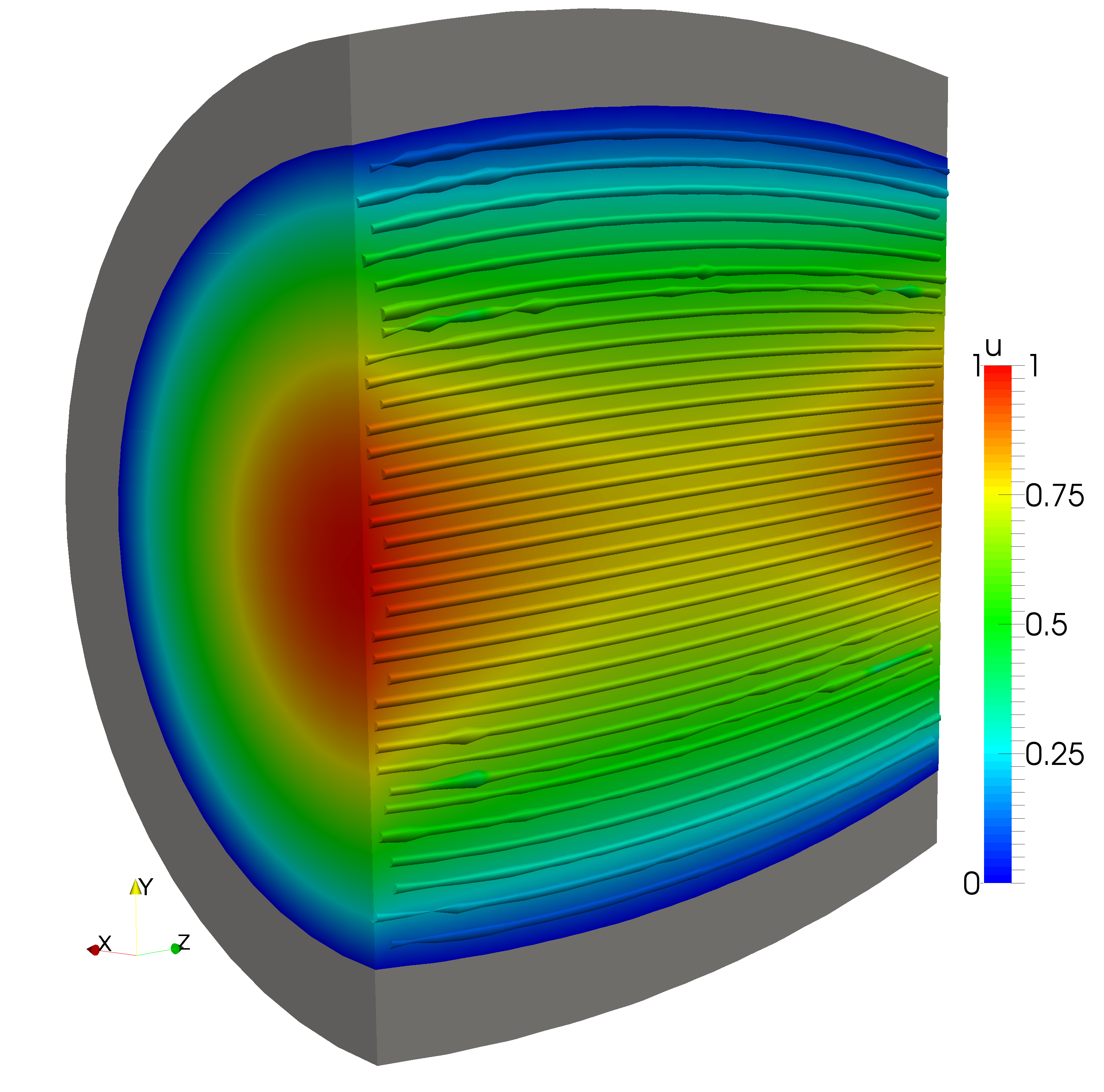}
  \hspace{1ex}
  \includegraphics[width=0.45\textwidth]{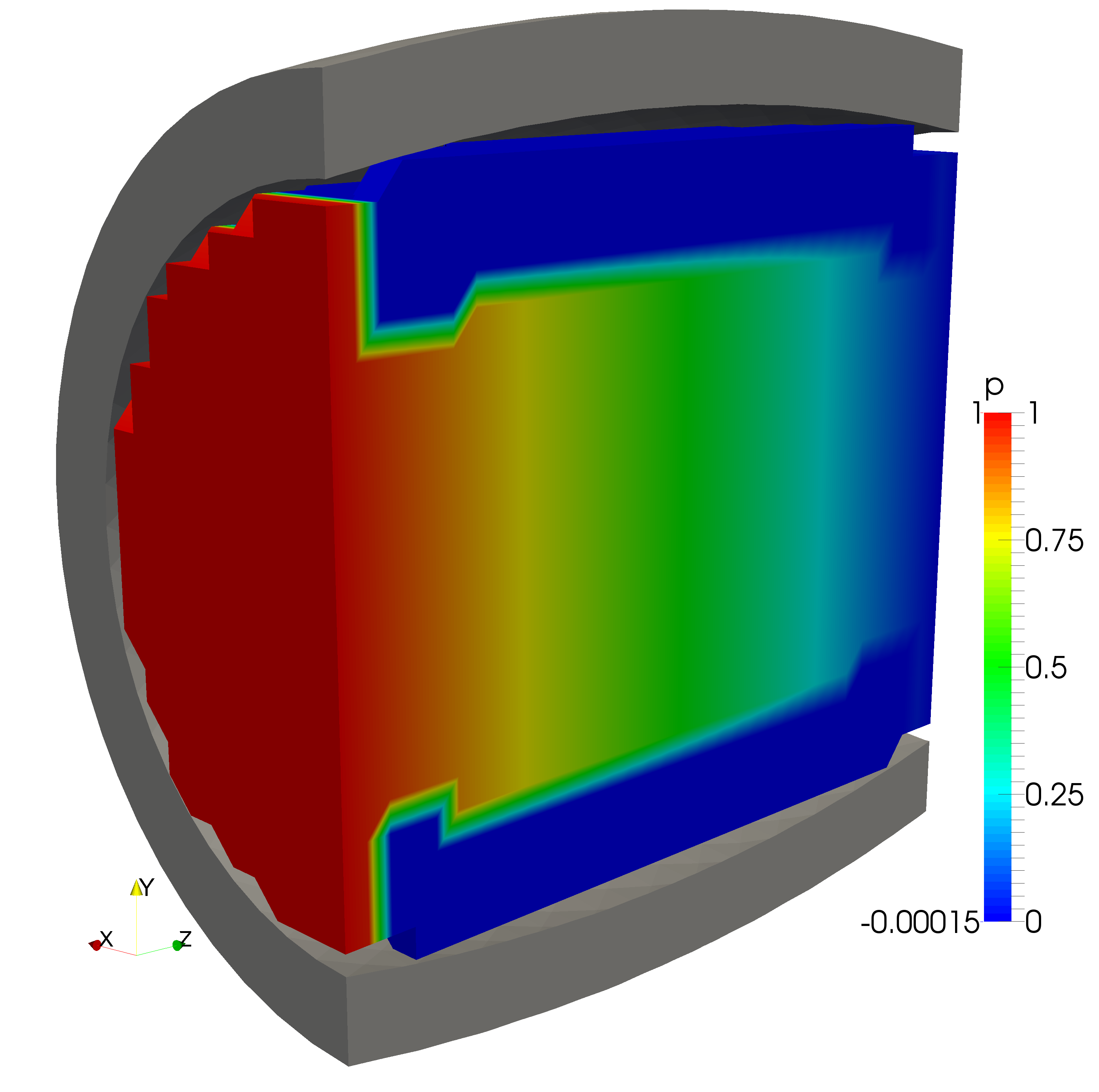}
  \includegraphics[width=0.45\textwidth]{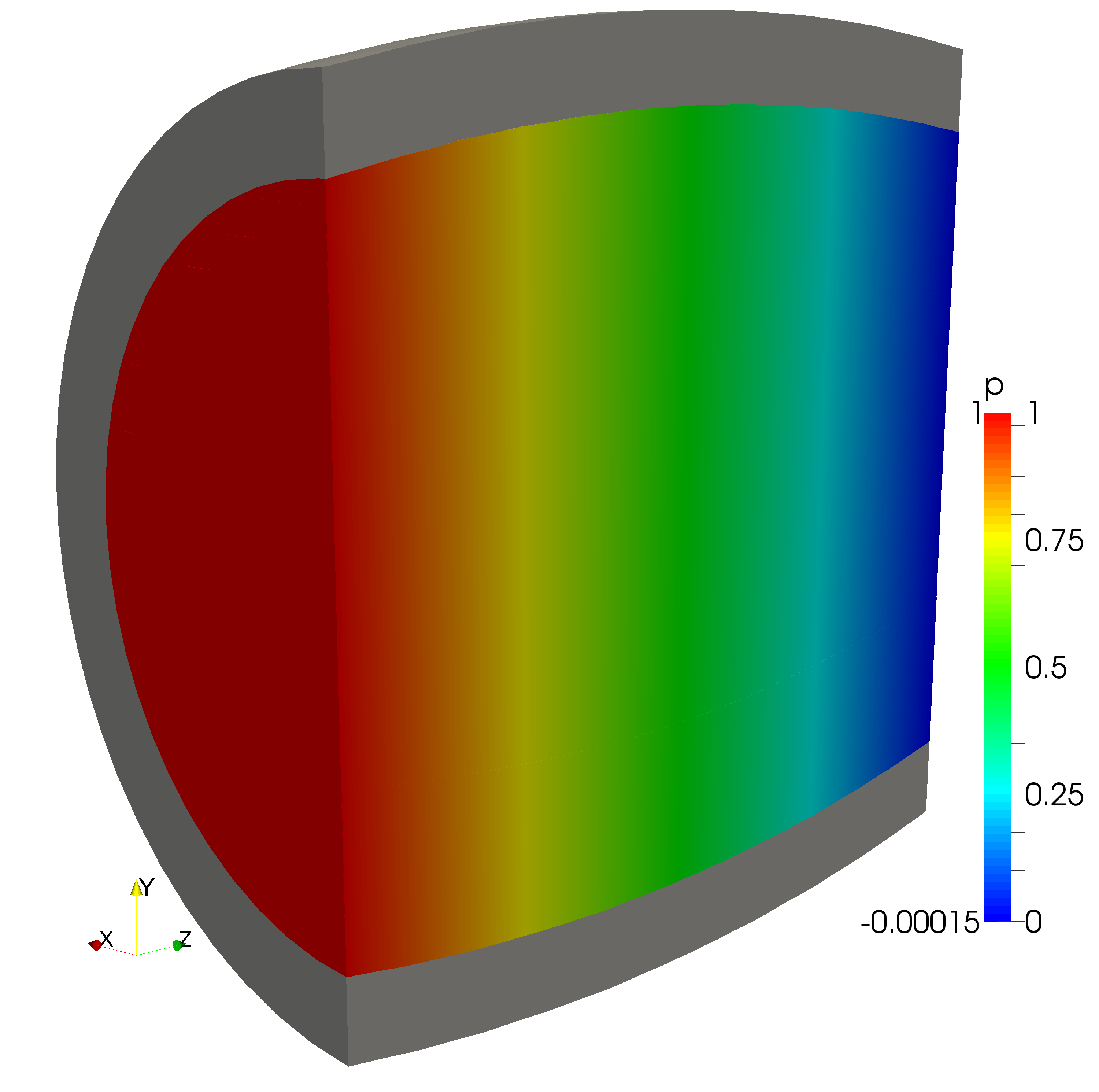}
  \caption{Computed velocity (top) and pressure (bottom) solution on
    the fixed fluid background mesh $\mesh^f_1$ (left) and entire
    overlapping fluid mesh $\{\mesh^f_1, \mesh^f_2\}$ (right) for the
    analytical FSI problem.}
  \label{fig:analytical-fsi-example-velocity-pressure}
\end{figure}

In the reference configuration, a discretization of the solid domain
$\refe{\Omega}^s$ and the fluid domain $\refe{\Omega}^f_2$ is provided
by two fitted and conforming meshes $\refe{\mesh}^s$ and
$\refe{\mesh}^f_2$, respectively, while the fluid domain
$\refe{\Omega}^f_1$ is represented by a structured Cartesian mesh
$\refe{\mesh}^f_1$ overlapped by the mesh $\refe{\mesh}^f_2$, see
Figure~\ref{fig:analytical-fsi-example-velocity-pressure}. The
numerical approximation of the fluid velocity, fluid pressure and
structure displacement are then computed on a sequence of $4$
overlapping meshes. The mesh sizes of the initial meshes
$\refe{\mesh}^f_1$, $\refe{\mesh}^f_2$, and $\refe{\mesh}^s$ are
$0.246$, $0.14$ and $0.212$, respectively and each of the subsequent
meshes is generated from the previous one by uniformly refining each
mesh. Based on the manufactured exact solution, the experimental order of
convergence (EOC) is then computed by
\begin{align*}
    \text{EOC}(k) = \dfrac{\log(E_{k-1}/E_{k})}{\log(2)}
\end{align*}
where $E_k$ denotes the error of the numerical solution computed at refinement level $k$.
The numerical experiment was conducted using $\nu^f = 0.001$
for the fluid viscosity and Lam\'e parameters given by
\begin{equation}
  \mu = E/(2 + 2\nu), \quad \lambda = E \cdot \nu /((1 + \nu)(1 - 2\nu))
  \label{eq:lame-constants}
\end{equation}
in $\Omega^s$ with $E = 10$ and $\nu = 0.3$.

\begin{figure}[htb]
  \label{fig:analytical-fsi-example-displacements}
  \centering
  \includegraphics[width=0.45\textwidth]{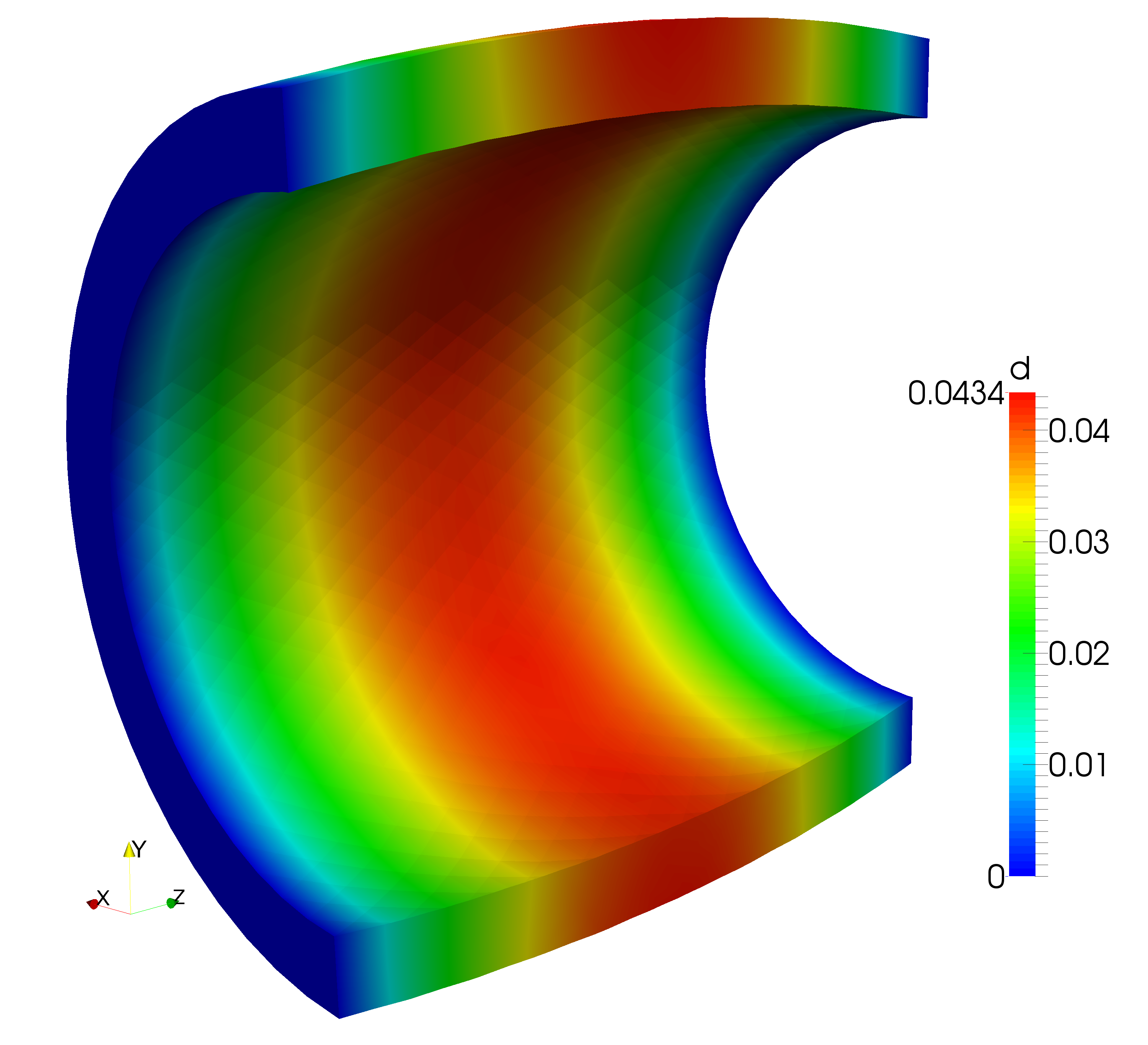}
  \includegraphics[width=0.45\textwidth]{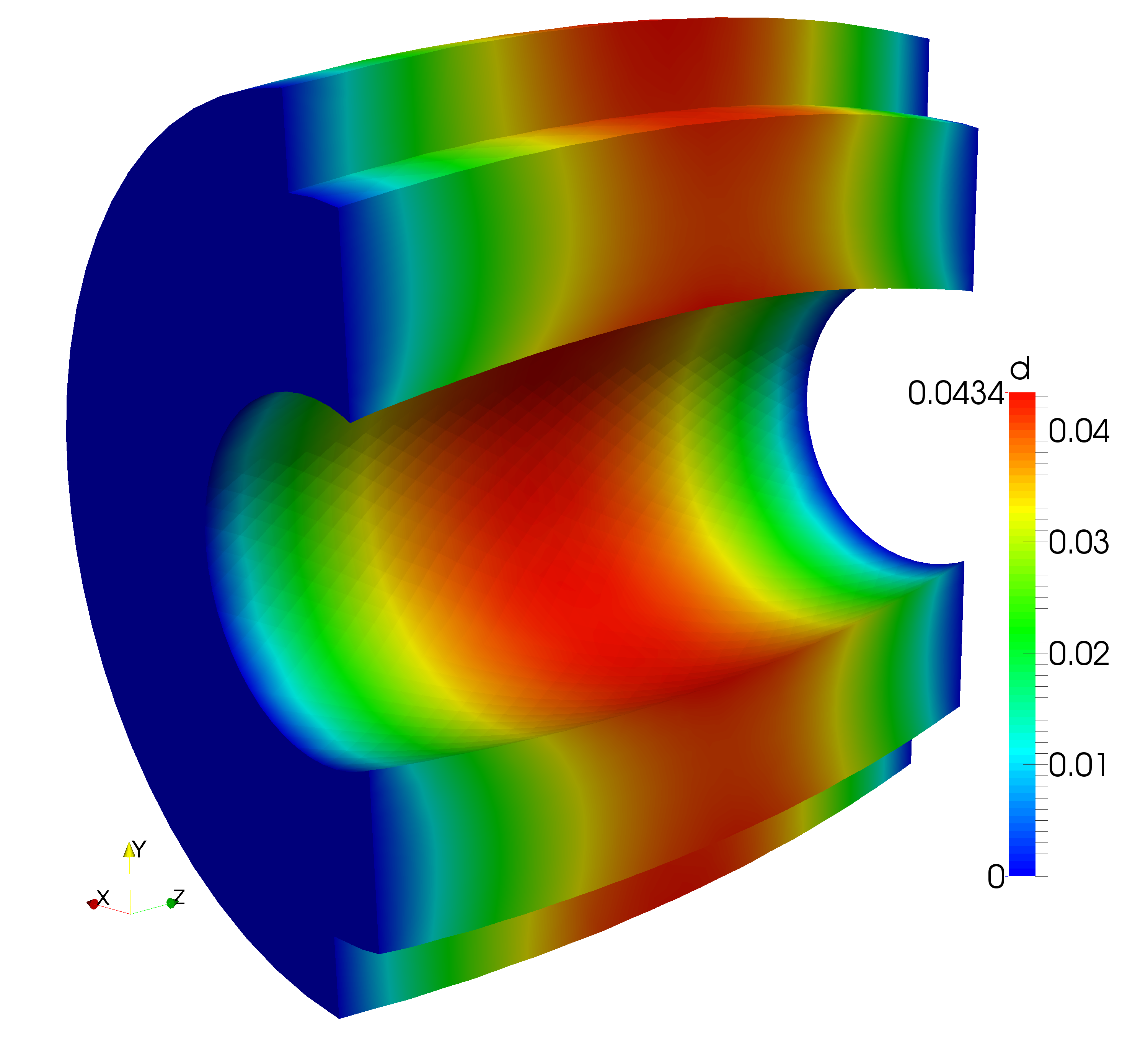}
  \caption{Displacements for the analytical FSI reference
    problem. Left: Structure displacement of the solid tube. Right:
    Displacement of the fluid mesh added.}
\end{figure}
For the penalty parameters in the stabilized overlapping mesh method
for the fluid problem, we pick $\gamma=10$ and $\delta = 0.5$.  Since
the overall computational time is dominated by the assembly and
solution of the fluid system, the displacement field is conveniently
solved using a direct solver, while the linear system arising from the
fluid problem is solved by applying a transpose-free quasi-minimal
residual solver with an algebraic multigrid preconditioner.

Using continuous piecewise linear functions for the approximation of
both the fluid velocity, fluid pressure and the structure
displacement, the theoretically predicted convergence rate for a
corresponding uncoupled problem is at least $1.0$ when measuring the
velocity and displacement error in the $H^1$-norm and the pressure
error in the $L^2$-norm. Note that it is common to observe a higher
experimental order of convergence of $\sim 1.5$ for the pressure
approximation when stabilized, equal-order interpolation elements are
used to discretize the flow problem. Assuming at most quadratic
convergence of the displacement solution in the $L^2$-norm, the
$L^2$-error will be reduced by approximately $0.5^{2\cdot3} \approx
0.016$ after 3 uniform mesh refinements.  To not pollute the overall
convergence rate by the iteration error, we therefore chose $tol =
0.001$ for the relative $L^2$-error between two consecutive
displacement solutions computed in the iteration loop. With the given
tolerance, the Dirichlet-Neumann iteration converged after $5$-$7$
iteration for each refinement level.  The resulting errors for the
sequence of refined meshes are summarized in
Table~\ref{tab:convergence-rates}. For the fluid velocity and fluid
pressure the observed convergence rates are in agreement with the
theoretical error decrease expected from an uncoupled problem.
For the solid displacement the observed convergence rates between
$1.46$-$1.9$ for the $H^1$-error are better than the
theoretically expected rate of $\sim 1$.

\begin{table}[htb]
  \centering
  \begin{tabular}{l|lr|lr|lr}
    \toprule
    Refinement &  $\| u^f_h - u^f \|_{1} $ 
               &  EOC &  $\| p^f_h - p^f \|_{0} $ &  EOC & $\| u^s_h -
    u^s \|_{1} $ &  EOC \\
    \midrule
    0 & $1.01188$ & -    & $3.61948e-03$  & - &    $3.87181e-04$ & -    \\
    1 & $0.51000$ & 0.99 & $1.55216e-03$ & 1.22 &  $1.40771e-04$ & 1.46 \\
    2 & $0.21912$ & 1.22 & $3.70746e-04$ & 2.06 &  $4.39062e-05$ & 1.68 \\
    3 & $0.12485$ & 0.81 & $1.29430e-04$ & 1.52 &  $1.17800e-05$ & 1.9  \\
    \bottomrule
  \end{tabular}
  \vspace{1ex}
  \caption{Convergence rates of the overlapping mesh finite element method for
    the analytical FSI problem.}
  \label{tab:convergence-rates}
\end{table}

\subsection{Flow around an elastic flap}
In the second numerical example, we consider a channel flow around an
elastic flap for different orientations of the flap with respect to
the channel geometry.  Here, the developed method and techniques can
play out their full strength as the overlapping mesh approach handles
large deformation within a single simulation easily.
As an additional benefit, our proposed scheme
allows to seamlessly reposition the flap for a series of numerical
experiments and thus has great potential for future applications
in design and optimization processes which involve fluid--structure interaction
problems in their forward simulation, see
for instance ~\cite{LundMollerJakobsen2003, EguzkitzaHouzeauxCalmetEtAl2014}.

Within the channel domain $\Omega = [0,L] \times [0,W] \times [0,H]$
with $L = 2.5$, $W = H = 0.41$, the bottom side of the flap of
dimensions $L^s = 0.06$, $W^s=0.2$, $H^s=0.24$ is centered around the
point $(L/2, W/2, 0)$.  In the first numerical experiment, the flap is
clamped on the boundary $[(L-L^s)/2,(L+L^s)/2],[(W-W^s)/2,(W+W^s)/2]
\times \{0\}$, while the flap is rotated $65^{\circ}$ around the
$z$-axis in a second experiment.  For the numerical experiment, we
assume that the flow can be described by the Stokes equations with
fluid viscosity $\nu^f = 0.001$, while the flap is modeled as an
hyperelastic material satisfying the St.~Venant--Kirchhoff
constitutive equation~\eqref{eq:venant-kirchhoff-material} with the
Lam\'e constants $\mu,\,\lambda$ defined by~\eqref{eq:lame-constants}
for $E^s = 15$ and $\nu^s = 0.3$.  Finally, we prescribe the inflow
profile $\bfu^f = (16\cdot0.45y(W-y)z(H-z), 0, 0)$ at the inlet $\{0\}
\times [0,W] \times [0,H]$, a ``do-nothing''
boundary condition given by $\nu
\partial_{\bfn} u - p\bfn = 0$ at the outlet $\{L\} \times [0,W]
\times [0,H]$ and a no-slip condition $\bfu = \bfzero$ elsewhere on
the boundary.

\begin{figure}[htb]
  \begin{center}
    \includegraphics[width=0.90\textwidth]{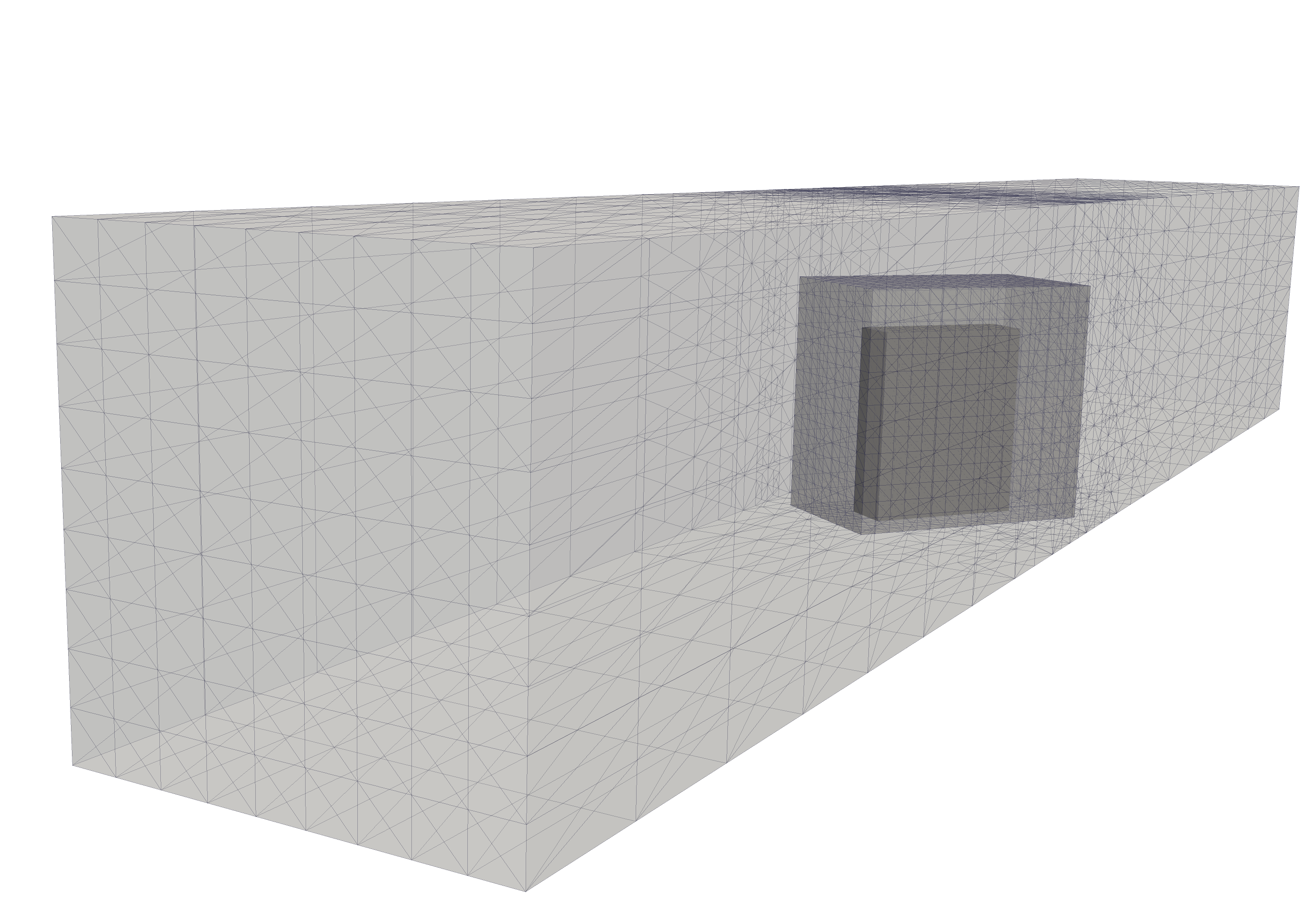}
    \caption{Background fluid mesh, structure mesh and its surrounding
      fluid mesh in the reference configuration.}
    \label{fig:meshes}
  \end{center}
\end{figure}

\begin{figure}
  \begin{center}
    \includegraphics[width=0.51\textwidth]{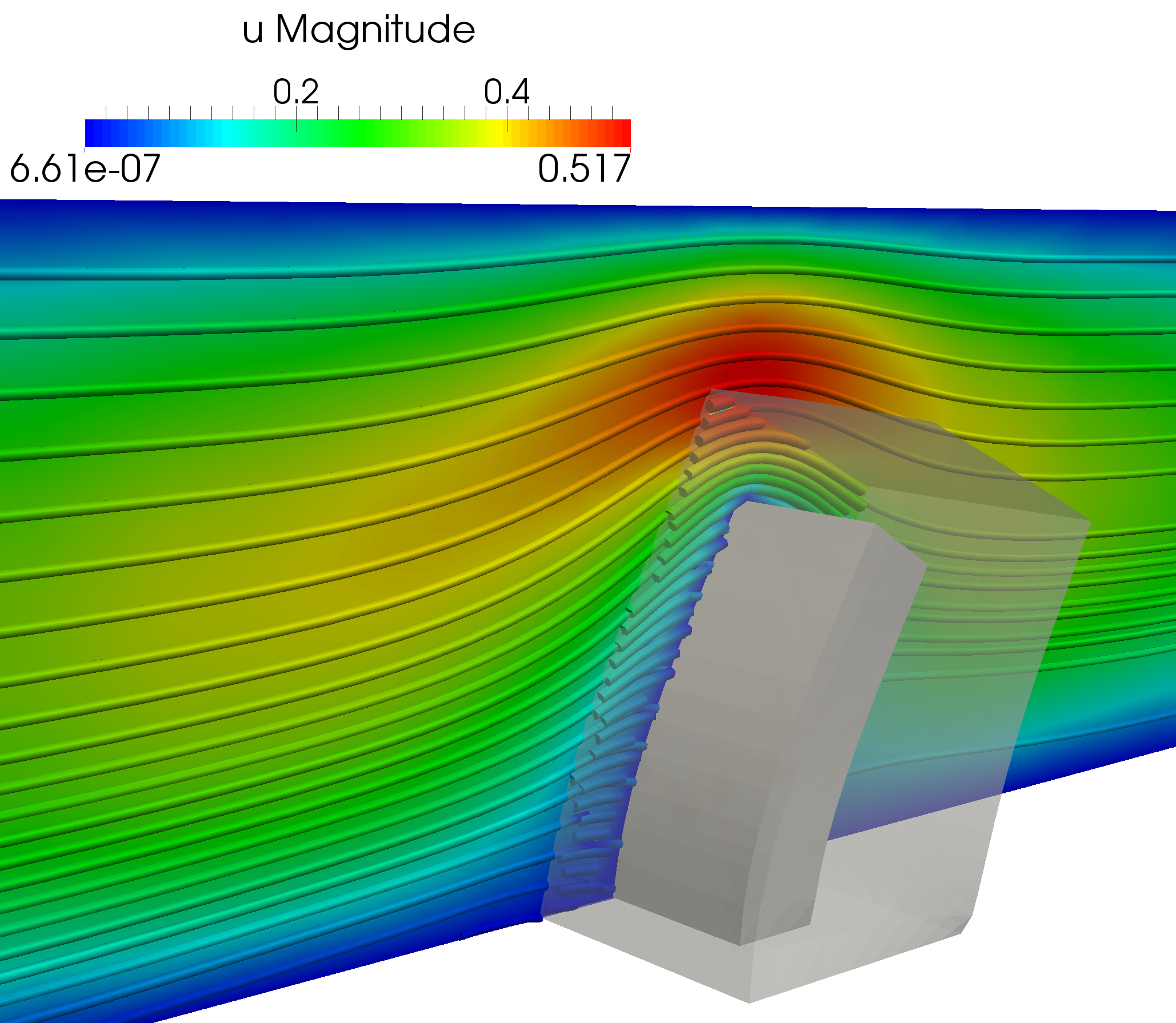}
    \includegraphics[width=0.48\textwidth]{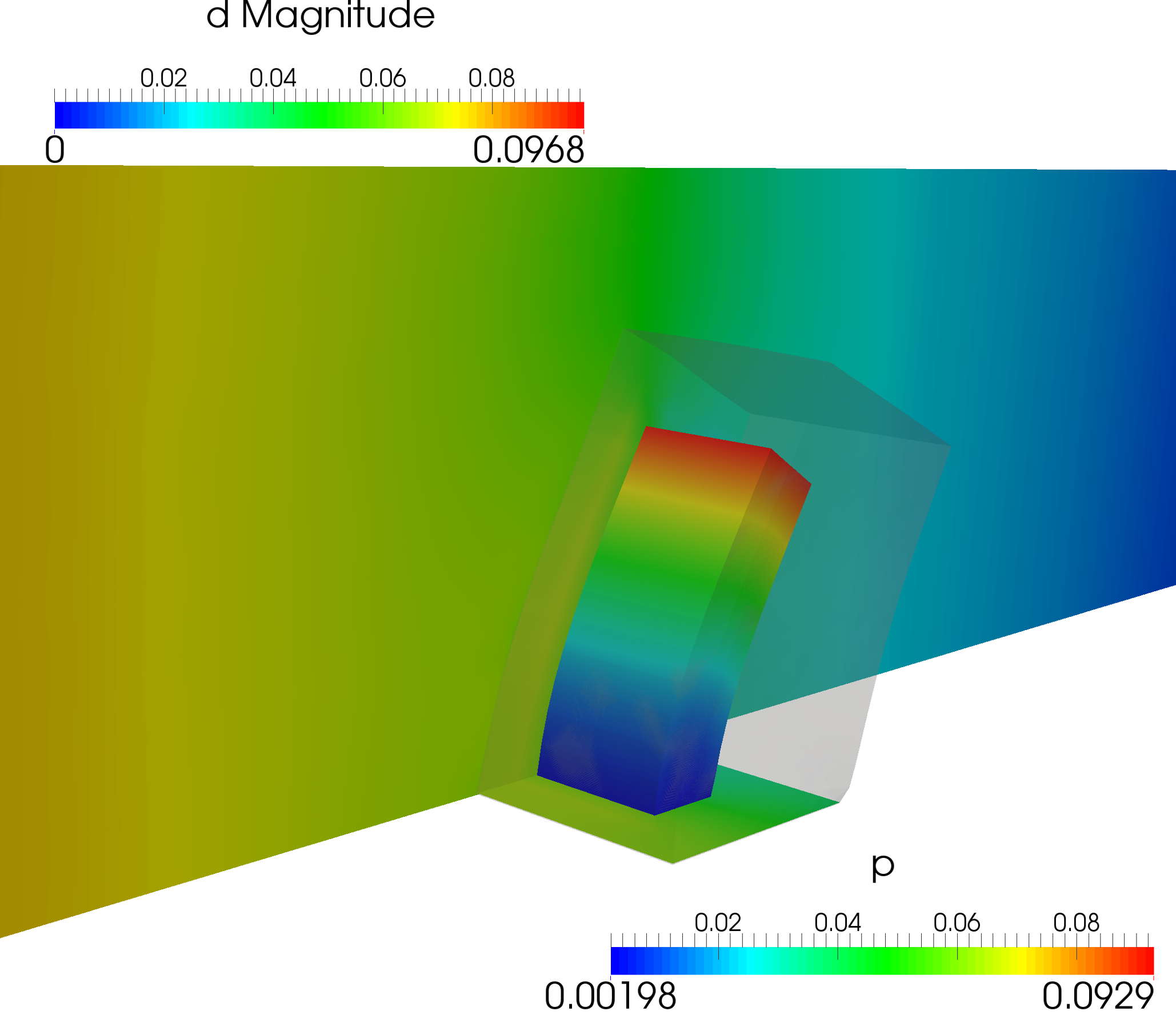}
    \\
    \includegraphics[width=0.52\textwidth]{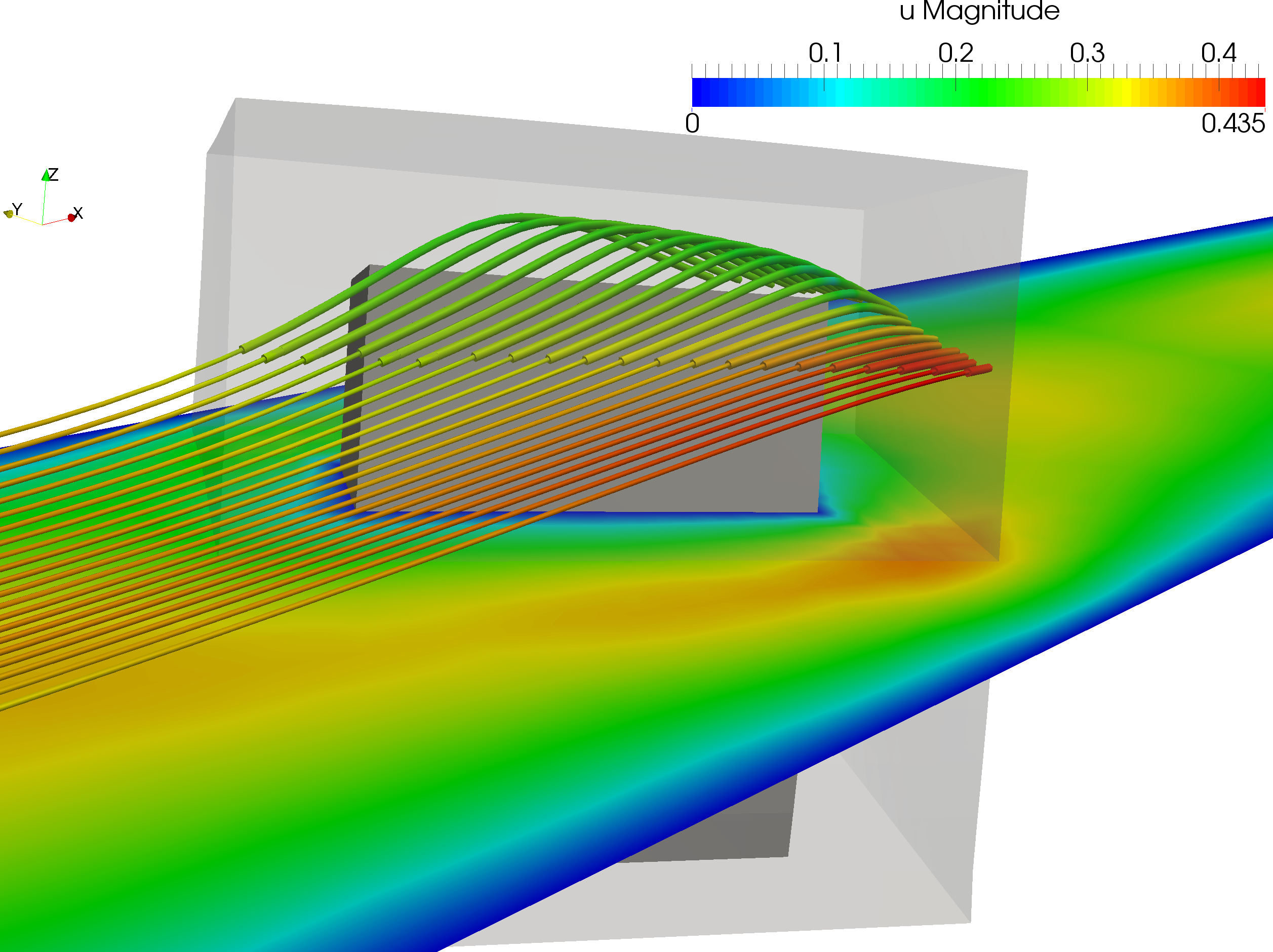}
    \includegraphics[width=0.47\textwidth]{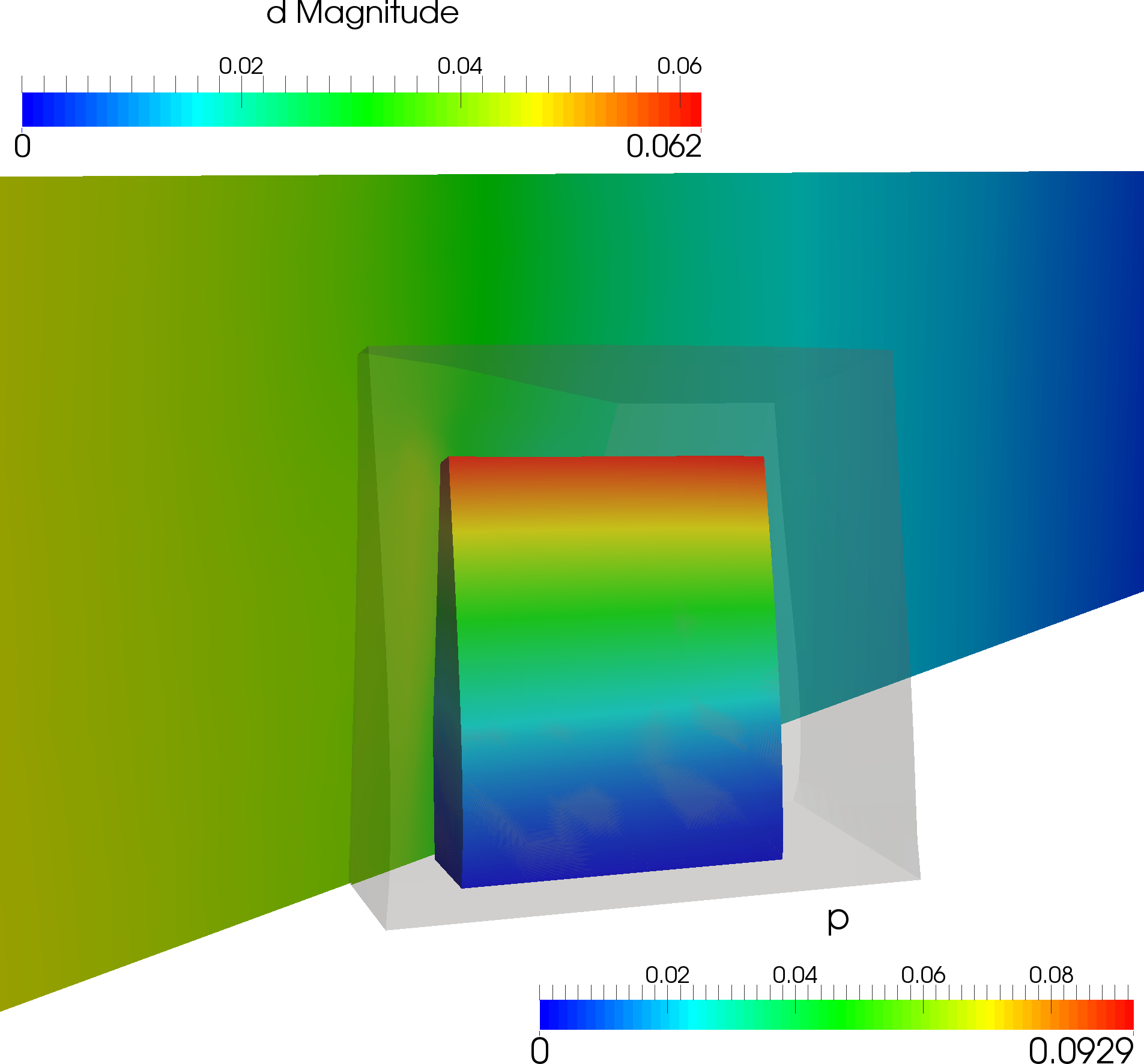}
    \caption{Flow around an elastic flap for two different
      flap orientations. Left:
      Magnitude and streamlines of the velocity approximation in an
      $x$-$z$ (top) and $x$-$y$ (bottom) cross-section.
      The transparent block around the
      gray-colored flap visualizes the fluid mesh $\mesh^f_2$
      surrounding the structure. The streamlines within the
      $\mesh^f_2$ are drawn slightly thicker to illustrate the smooth
      transition of the velocity approximation from the outer to the
      inner fluid domain. Right: Pressure distribution and magnitude
    of the structure displacement.}
    \label{fig:flap-channel-solutions-angle-0} \end{center}
\end{figure}

The numerical results for aligned and rotated flaps are shown in
Figure~\ref{fig:flap-channel-solutions-angle-0}.
We especially note the smooth transition of the velocity and pressure
solutions from fluid background $\mesh^f_1$ to the solid-surrounding
fluid mesh $\mesh^f_2$; the interface is not visible. The meshes used
for simulation of the rotated flap are shown in
Figure~\ref{fig:meshes}.

\section{Conclusions}
We presented a Nitsche-based cut and composite mesh method for fluid--structure
interaction problems. The method utilizes a Nitsche type coupling
between two fluid meshes: one fixed background mesh and one moving
overlapping fluids mesh which is fitted to the boundary of a
hyperelastic object and deforms with the object. The fluid--fluid
coupling is monolithic in the sense that it manufactures a coupled
system involving both the underlying and overlapping degrees of
freedom. In previous work, \cite{MassingLarsonLoggEtAl2013}, we have shown that the
coupling is stable and that the solution has optimal order convergence
for a stationary model problem.

To solve for the steady state solution of a fluid--structure
interaction problem with large elastic deformations, we consider a
fixed-point iteration where we solve for the fluid, compute a boundary
traction for the solid, solve for the solid, solve for the mesh motion
of the overlapping fluid mesh, and finally update the geometry. This
involves computing new intersections between underlying and
overlapping meshes.  Employing a provably stable overlapping mesh
method for fluid-fluid coupling, the proposed scheme for the
fluid--structure problem is guaranteed to be robust and insensitive to
the overlap configuration.

We verified the expected convergence rates for a model problem with a
manufactured solution and demonstrated the flexibility of our approach
by computing the steady state solution for an elastic flap in a
channel at two different orientations. It should be noted that the
overlapping mesh method allows the flap to be repositioned in the
channel without requiring the generation of a single \emph{conforming}
fluid mesh for each configuration. Only an element-wise, local
representation of the cut cells near the interface together with some
appropriate quadrature schemes are required, see for
instance~\cite{Massing2012a}.

Future work involves extending our method to fully time-dependent flow
governed by the incompressible Navier-Stokes equations. We note that
the nonlinear convection term can be handled in our setting using a
discontinuous Galerkin coupling with up-winding and that, from a
computational point of view, taking a time step is closely related to
taking one step in our fixed-point iteration algorithm. Another area
of interest is the direct coupling between fluids and solids.

\section*{Acknowledgments}

This work is supported by an Outstanding Young Investigator grant from
the Research Council of Norway, NFR 180450. This work is also
supported by a Center of Excellence grant from the Research Council of
Norway to the Center for Biomedical Computing at Simula Research
Laboratory. The authors would like thank the anonymous referee for the
valuable comments and suggestions which helped to improve the
presentation of this work.

\bibliographystyle{plainnat}
\bibliography{bibliography}

\end{document}